\documentclass[final,leqno]{siamltex}

\usepackage[ruled]{algorithm2e}

\usepackage{subfigure,upgreek,dsfont,braket,amssymb,amsmath,booktabs}
\usepackage[small,bf]{caption}
\usepackage{dsfont}
\usepackage{float,url}
\usepackage{graphicx}
\usepackage{color,float}
\usepackage{epsfig}

{\makeatletter \@addtoreset{equation}{section}}

\newtheorem{remark}[theorem]{Remark}

\newcommand{\R}{ \mathbb{R} }

\newcommand{\norm}[2]{\left\lVert#1\right\rVert_{#2}}

\newcommand{\N}{\mathds{N}}

\newcommand{\cU}{\mathcal{U}}

\newcommand{\cI}{\mathcal{I}}

\newcommand{\sgn}{\operatorname{sgn}}

\newcommand{\D}{\partial}

{\begin{minipage}{0.9\textwidth}\hrule height.08em\kern4pt\begin{algorithmic}[1]}%
  {\end{algorithmic}\kern4pt\hrule height.08em\relax\end{minipage}}

\def\e{{\epsilon}}

\def\norm#1{\|#1\|}

\title{Local minimization algorithms for dynamic programming equations}

\author{Dante Kalise\thanks{Johann Radon Institute for Computational and Applied Mathematics
(RICAM), Austrian Academy of Sciences, Altenberger Stra\ss{}e 69, 4040 Linz, Austria ({\tt
dante.kalise@oeaw.ac.at}).}
\and Axel Kr\"oner\thanks{INRIA Saclay and CMAP, \'Ecole Polytechnique, Route de Saclay, 91128
   Palaiseau Cedex, France, ({\tt
 axel.kroener@inria.fr})}.
\and Karl Kunisch\thanks{University of Graz, Institute of Mathematics and Scientific
Computing, Heinrichstr. 36, A-8010 Graz, Austria and Johann Radon Institute for Computational and Applied Mathematics
(RICAM), Austrian Academy of Sciences, Altenberger Stra\ss{}e 69, 4040 Linz, Austria ({\tt
karl.kunisch@uni-graz.at}).
}}

\begin{document}

\maketitle

\begin{abstract}
The numerical realization of the dynamic programming principle for
continuous-time optimal control leads to nonlinear Hamilton-Jacobi-Bellman equations which require the minimization of a nonlinear mapping over the set of admissible controls. This minimization is often performed by comparison over a finite number of elements of the control set. In this paper we demonstrate the importance of an accurate realization of these minimization problems and propose algorithms by which this can be achieved effectively. The considered class of equations includes nonsmooth control problems with $\ell_1$-penalization which lead to sparse controls.
\end{abstract}



\pagestyle{myheadings}
\thispagestyle{plain}
\markboth{Dante Kalise, Axel Kr\"oner, Karl Kunisch}{ }


%

\section{Introduction}

Since its introduction by Bellman in the 50's, dynamic programming has become a fundamental tool in the design of optimal control strategies for dynamical systems. It characterizes the value function of the corresponding optimal control problem in terms of functional relations, the so-called Bellman and Hamilton-Jacobi-Bellman (henceforth HJB) equations. We begin by briefly recalling this setting in the context of infinite horizon optimal control.

\noindent We make the following assumptions. We equip $\R^n$ for $n\in \N$ with the
Euclidean norm $\norm{\cdot}_2$.
Furthermore, let 
\begin{equation}
\begin{aligned}
&f\colon \R^d\times \R^m \rightarrow \R, &  |f(x,u)-f(y,u)| & \le  \omega_R \norm{x-y}_2,\\
&l\colon \R^d \times \R^m \rightarrow \R,& |l(x,u)-l(y,u)|& \le  \omega_R \norm{x-y}_2,
\end{aligned}
\end{equation}
for $x, y\in\R^ d$ with $\|x-y\|_2\leq R$ and modulus $\omega_R\colon [0,\infty) \rightarrow [0,\infty)$ of polynomial growth 
satisfying $\lim_{r \rightarrow 0^+} \omega_R(r)=0$ (cf. Ishii
\cite{Ishii:1984}), $d,m\in \N$.
Let the dynamics be given by
\begin{eqnarray}\label{eq:dynconaut}
	\left\{\begin{array}{l}
	\dot{y}(t)=  f(y(t),u(t)),\cr
	y(0)=  x
	 \end{array}\right.
\end{eqnarray}
for $t>0$, where  $x\in \R^d$ and $u\in\cU\equiv\{ u\colon \R_+\rightarrow U\,\text{measurable}\}$,
$U\subset \R^m$ compact. We introduce the following cost functional $J:\cU\rightarrow\R $ 
\begin{eqnarray*}\label{eq:joi}
	J(u)=\int_0^{\infty} l(y(s),u(s))e^{-\lambda s}ds\,,\quad\lambda>0\,,
\end{eqnarray*}
where $y$ is the solution of \eqref{eq:dynconaut} depending on $x$ and $u$.
By the application of the dynamic programming principle, the value function 
\[v(x)\equiv \inf_{u\in\cU} J(u)\]
is characterized as the viscosity solution \cite[Chapter 3]{BardiCapuzzoDolcetta:2008} of the HJB equation
\[
\lambda v(x)+\sup_{u\in U}\{-f(x,u)\cdot \nabla v(x)-l(x,u)\}=0,\quad x\in \R^d.
\]
There exists an extensive literature concerning the construction of numerical schemes for static HJB equations. The spectrum of numerical techniques includes ordered upwind methods \cite{vlad,mitchell}, high-order schemes \cite{shu}, domain decomposition  techniques \cite{cacace} and geometric approaches \cite{botkin}, among many others (we refer to  \cite[Chapter 5, p.145]{FalconeFerretti:2014} for a review of classical approximation methods). In this paper we follow a semi-Lagrangian approach \cite{falconefirst}, which is broadly used to approximate HJB equations arising in optimal control problems, see,e.g., \cite{FalconeFerretti:2014}. We illustrate the basic steps in the formulation of a semi-Lagrangian scheme for our model problem.

The construction of a first-order semi-Lagrangian scheme begins by considering an Euler discretization
of the system dynamics with time step $h>0$
\begin{equation*}
	\left\{
        \begin{aligned}
	y^{n+1}&=y^n+hf(y^n,u^n),\\
	y^0&=  x,
      \end{aligned}
      \right.
\end{equation*}
for $n\in \N^0$, $x\in \R^d$, and controls $u^n\in U$. Then, the application of the dynamic programming principle on the discrete-time dynamics leads to the Bellman equation
\begin{align*}
 v(x)=\min_{u\in U}\{(1-\lambda h)v(x+hf(x,u))+hl(x,u)\},\, \quad x\in \R^d.
\end{align*}

\noindent To discretize this equation in space we introduce a bounded domain $\Omega=[a,b]^d\subset \R^d$, $a,b\in \R$, where we define a regular quadrangular mesh with $N$ nodes and mesh parameter $k$. We denote the set of nodes by $\Omega_k\subset \Omega$. 
Let the discrete value function be defined in all grid points, $V:=\{v(x)\}_{x \in \Omega_k}$. 
However, note that $x+f(x,u)$ for $x \in \Omega_k$ is not necessarily a grid point, 
and therefore the value function has to be evaluated by interpolation which is chosen as a linear one
here. The interpolant $I[\cdot](x)$ is
defined on the basis of the dataset $V$. The resulting fully discrete scheme
then reads
\begin{equation}\label{SLih}
 V(x)=\min_{u\in U}\{(1-\lambda h)I[V](x+hf(x,u))+hl(x,u)\}=:G(V), \quad x \in  \Omega_k,
\end{equation}
which can be solved by a fixed point iteration starting from an initial guess $V ^0$ by
\[
V^{i+1}=G(V^i)\,.
\]
 An alternative to the fixed point approach is the use of Howard's algorithm or iteration in the policy (control) space \cite{boka,AllaFalconeKalise:2013},  which is faster but also utilizes the parametric minimization of the Hamiltonian.

\noindent Finally, once the value function is computed this allows to derive a feedback control for a given
state $x \in \R^d$ by
\[
u(x)\in \underset{u\in U}{\operatorname{argmin}}\{(1-\lambda h)I[V](x+hf(x,u))+hl(x,u)\}.
\]

A characteristic feature of the class of HJB equations arising in optimal control problem is its
nonlinear Hamiltonian, 
\[
\sup_{u\in U}\{-f(x,u)\cdot \nabla v(x)-l(x,u)\}\,,
\]
which requires a parametric maximization (minimization) over the control set $U$. For its discrete analogue $G(V)$, a common practice in the literature is to compute the minimization by comparison, i.e., by evaluating the expression in a finite set of elements of $U$ (see for instance \cite{AllaFalconeKalise:2013,falconegame,KroenerKunischZidani:2013} and references therein). In contrast to the comparison approach, the contribution of this paper is to demonstrate that an accurate realization of the min-operation on the right hand side of \eqref{SLih} can have an important impact on the optimal controls that are determined on the basis of the dynamic programming principle. In this respect, the reader can take a preview to Figure \ref{test1},  where differences between optimal control fields obtained with different minimization routines can be appreciated. Previous works concerning the construction of minimization routines for this problem date back to \cite{CarliniFalconeFerretti:2004}, where Brent's algorithm is proposed to solve high dimensional Hamilton-Jacobi -Bellman equations and to \cite{CristianiFalcone:2007}, where the authors consider a fast  semi-Lagrangian algorithm for front propagation problems.  In this latter reference, the authors 
determine the minimizer of a specific Hamiltonian by means of an explicit formula. 
Moreover, for local optimization strategies in dynamic programming we refer to  \cite{Jarczyk:2005} for Brent's algorithm  and to \cite{Unrath:2006} for a Bundle Newton method. 

In this article, a first-order primal-dual method (also known as
Chambolle-Pock algorithm \cite{ChambollePock:2011}) and a semismooth Newton method \cite{ItoKunisch:2008, Ulbrich:2011} are proposed within the semi-Lagrangian scheme for the evaluation of the right hand side in \eqref{SLih}. In contrast to the minimization by the comparison approach, the
proposed algorithms leads to more accurate solutions  for the same CPU time. Since we preserve the continuous nature of the control set, in some specific settings it is also possible to derive convergence results for our minimization strategies. Furthermore, it provides a solid framework to address challenging issues, such as nonsmooth optimal control
problems with $\ell_1$ control penalizations in the cost functional.

%
%

The paper is organized as follows. In Section \ref{sec:interpolation_operator} we 
begin by recasting the discretized Hamiltonian as a minimization problem explicitly depending on the control $u$. In Section \ref{sec:alg} we introduce and adapt the Chambolle-Pock and semismooth
Newton methods for the problem under consideration, and in Section \ref{sec:examples} we present numerical
examples assessing the accuracy and performance of the proposed schemes.

\section{Explicitly control-dependent Hamiltonians}\label{sec:interpolation_operator}
In this section we study the numerical evaluation of the discretized nonlinear Hamiltonian 
\begin{equation}\label{numham}
\min_{u\in U}\{(1-\lambda h)I[V](x+hf(x,u))+hl(x,u)\}\,.
\end{equation}
This is a non -- standard minimization problem requiring the evaluation of the nonlinear mapping $u\rightarrow I[V](x+hf(x,u))$ which depends on the discrete dataset $V$ and the system dynamics $f(x,u)$. Therefore, a first step towards the construction of local minimization strategies is to recast \eqref{numham} by assuming specific structures for $I[V], f(x,u),$ and $l(x,u)$, leading to explicit piecewise linear or quadratic optimization problems on $u$. This section is split between the treatment of the interpolation $I[V](x+hf(x,u))$, and the evaluation of $l(x,u)$. For the sake of simplicity we restrict the presentation to the two dimensional case $d=2$, although the presented framework can be generalized to
higher dimensions.

\subsection{The interpolation operator and local subdivisions of the control space}

\begin{figure}[!ht]
  \centering
  \epsfig{file=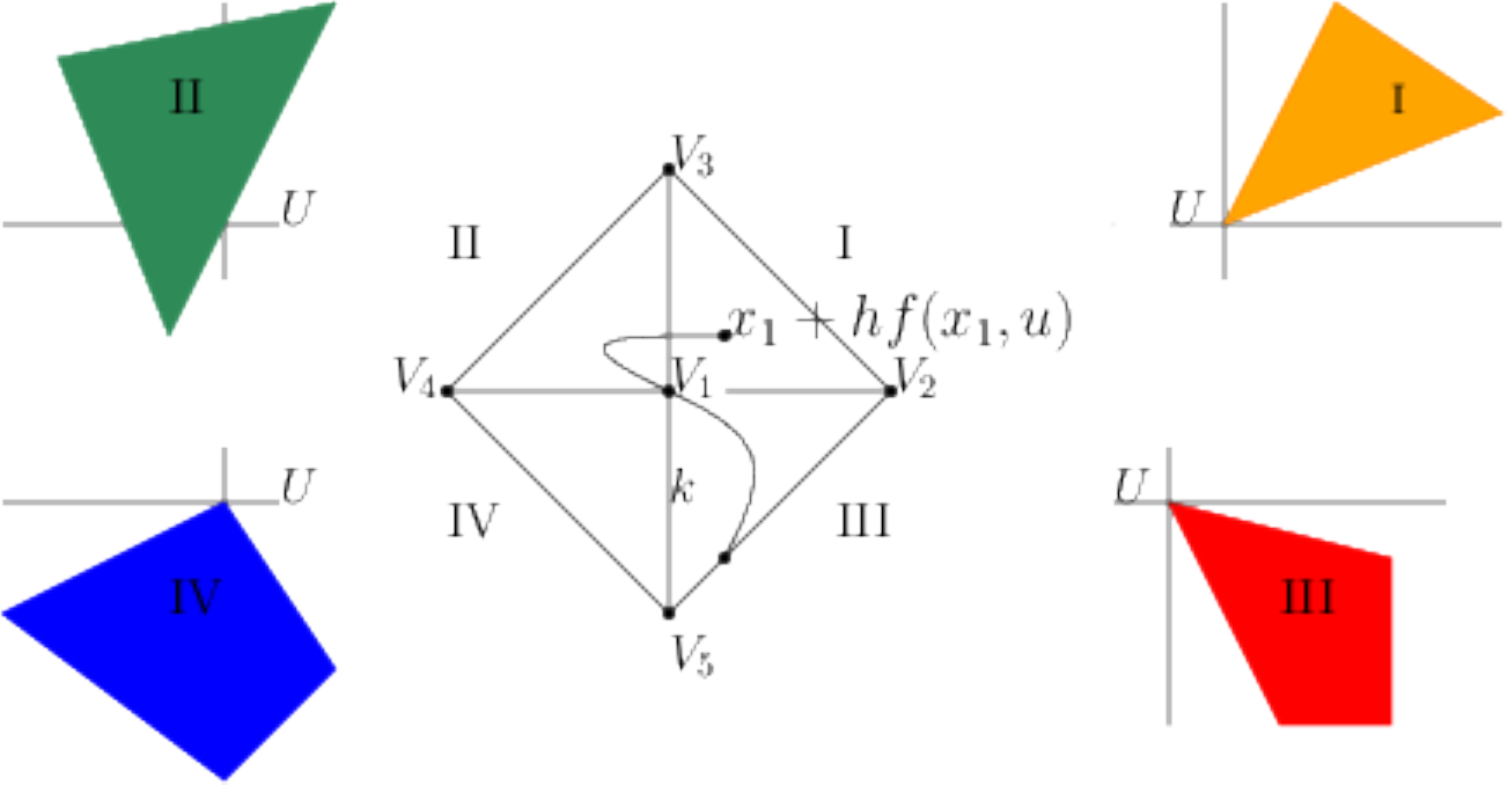,width=.7\linewidth,clip=}
  \caption{Arrival points and related control sets.}\label{figinterp}
\end{figure}

To analyze the interpolation operator  we consider for every node $x=(x_1,x_2) \in \Omega_k$ the patch of four
triangles defined by the neighboring nodes, cf. Fig.~\ref{figinterp}. Assume that the arrival point $x+hf(x,u)$ for $x \in \Omega_k$, $u\in
U$, is located in a triangle defined by the 
points $x=x^1$, $x^2$,  and $x^3$ in $\R^2$ with associated values $V_1$, $V_2$, and $V_3$ as indicated
in Fig.~\ref{figinterp}. The linear
interpolation formula then reads for a mesh point $x\in \Omega_k$
\begin{equation}\label{bil}
I[V](x)     =  c x_1 + d x_2 +e\,
\end{equation}
with
\begin{equation}
\left(\begin{array}{ccc}
x_1^1&x_2^1&1\\
x_1^2&x_2^2&1\\
x_1^3&x_2^3&1\\
\end{array}\right)
\left(\begin{array}{c}
c\\
d\\
e\end{array}\right)
=
\left(\begin{array}{c}
V_1\\
V_2\\
V_3\end{array}\right).
\end{equation}

\begin{remark}
{ \em Note that in the considered case of first order approximations the passage from the interpolation
on one triangular sector to another one is continuous. For the first order interpolation schemes in a point $x\in \Omega_k$, one might also think of considering one interpolant for a macro-cell defined by a larger set of neighboring nodes. This could be
computed by using $2^{d}$ suitable points. However, in practice this is similar to consider a grid
 of size $2h$, which leads to less accuracy, which is, in particular, in higher dimensions a
crucial issue. Alternatively one could consider one interpolant for the macro-cell which is defined
piecewise over all the triangulars. However, this piecewise smooth interpolation is not convenient for
numerical optimization. As a further alternative one can resort to higher order interpolation
schemes.
}$\square$\end{remark}

 According to \eqref{numham}, the interpolation operator $I[V](x)$ is evaluated at the arrival point $x+hf(x,u)$. If this point is sufficiently close to $x\in \Omega_k$,   it is in general possible to know a priori in which of the four triangular sectors it will be located. Our goal is to establish a consistent procedure to locally divide the control space in sets that generate arrival point related to a single triangular sector. For instance, in the case of the eikonal equation there is a simple correspondence between the computed triangulars in control and state spaces, see Section \ref{sec:eikonal}.
In the following 
we restrict ourselves to dynamics of the form
  \begin{align}\label{form_dynamics}
f(x,u)=g(x)+Bu,
\end{align}
for $ B \in \R^{d\times m}$, i.e. nonlinear in the state, and linear in the control. 
 We assume that $\bar h >0$ is chosen such that
 \begin{align}\label{mesh_paramter_cond}
\overline{h} \, \underset{x \in \Omega, u \in U}{\text{sup}} \|g(x) + Bu\|_2 \le
\frac{\sqrt{2}}{2}  k.
\end{align}
Let the rows of $B$ be denoted by $\{b_i\}_{i=1}^d$.
Further, let $x_D$ denote a grid point in $\Omega_k$, let $\mathcal{I} \subset \{1,\dots,d\}$ be an index set
with complement $\mathcal{I}^\subset$, and define a triangle $Q_{\mathcal{I}} \subset \R^d$ associated to $x_D$
for the interpolation and $\mathcal{I}$ by
\[
Q_{\mathcal{I}} =\{\tilde{x} = x_D + x \, | \, \norm{x}_1 \le k,\, x_i \ge 0 \, \text{ for }\, i \in
\mathcal{I},\,x_i \le 0 \, \text{ for }\, i \in \mathcal{I}^\subset \}.
\]
Then
\[
\mathcal{U}_\mathcal{I} = \{u \in U\, | \; g(x_D)_i + b_iu \ge 0 \, \text{ for }\, i \in
\mathcal{I},\, g (x_D)_i + b_i u \le 0 \, \text{ for }\, i \in \mathcal{I}^\subset\}
\]
has the property that $x_A = x_D + h(g(x_D) + Bu) \in Q_{\cI}$ for $u \in \mathcal{U}_\mathcal{I}$, and $0
\le h \le \overline{h}$. In fact, $x = x_A - x_D = h(g(x_D) + Bu)$, and we obtain with
\eqref{mesh_paramter_cond} that $\norm{x}_1 \le k,\, x_i \ge
0$ for $i \in \mathcal{I}$ and $x_i \le 0$ for $i \in \mathcal{I}^\subset$. Note also that 
\[
U =
\underset{\mathcal{I} \in \mathcal{W}} \bigcup \mathcal{U}_\mathcal{I},
\]
where $\mathcal{W}$ is the set of all subsets of index sets in $\{1,\dots,d\}.$ The associated interpolation operator on $Q_{\cI}$ is denoted by $I_{\cI}$. Note that evaluation of $I_{\cI}[V](x+f(x,u))$ for $x\in Q_{\cI}$ and $u\in\mathcal{U}_\mathcal{I}$ leads to a linear dependence on $u$ of the form
\begin{equation}\label{interpterm}
c_\cI u_1+d_{\cI} u_2+e_{\cI}\,,
\end{equation}
where the coefficients $c_\cI$, $d_\cI$ and $\e_{\cI}$ are uniquely determined by the vertices of $Q_\cI$.

\subsection{Evaluation of the cost term}
Having approximated the interpolation term of the Hamiltonian, what is left is to provide an
approximation of the running cost $l(x,u)$. 
If this term is defined in a pointwise manner, this imposes no additional difficulty. For instance,
in some of the examples we shall utilize
\begin{align}\label{l_example}
l(x,u)=\|x\|_2^2+\frac{\gamma}{2}\|u\|_2^2=x_1^2+x_2^2+\frac{\gamma}{2}(u_1^2+u_2^2)\,\quad x\in
\Omega_k,\quad u\in U.
\end{align}
This introduces a constant and a quadratic term in $u$ to be added to the above presented
expressions for the interpolant \eqref{interpterm}. For $l$ given as in \eqref{l_example} we will
consider in the scheme \eqref{SLih} the explicit
dependence on $u$. For a given node $x\in \Omega_k$ we have
\begin{align}
  [V]_{x}&=\underset{\cI\in \mathcal W}{\min}\{H_{\cI}(x,V,u,\cI)\}\,,\label{nh1}
\end{align}
where $[V]_x$ denotes the value of the discrete value function at node $x$. 
Further, for each set $\cI$ we have
\begin{align}
  H_{\cI}(x,V,u,\cI)&=\underset{u\in U_{\cI}}{\min}\{\beta I_{\mathcal I}[V](x+hf(x,u))+hl(x,u)\}\,,\label{nh2}\\
 &=\underset{u\in
 U_{\cI}}{\min}\left\{\beta(c_\cI u_1+d_{\cI} u_2+e_{\cI})+h\left(x^2+y^2+\frac{\gamma}{2}(u_1^2+u_2^2)\right)\right\}\,,\nonumber\\\label{H}
 &=\underset{u\in U_{\cI}}{\min}\{\frac{\tilde a_{\cI}}{2}u_1^2+\frac{\tilde b_{\cI}}{2}u_2^2+\tilde
 c_{\cI}u_1+\tilde d_{\cI}u_2+\tilde e_{\cI}\}\,,
\end{align}
with
\begin{align*}
  &\beta =1-\lambda h\,, \tilde a_{\cI}=h\gamma\,,
  \tilde b_{\cI}= h\gamma\,,\;  \tilde c_{\cI}=\beta c_{\cI}\,,\\
  &\tilde d_{\cI}=\beta
 d_i\,, \tilde e_{\cI}=\beta e_{\cI}+ h (x^2+y^2)\,.
\end{align*}



\begin{remark}
  {\em For minimum time optimal control problems, the resulting semi-Lagrangian scheme for points $x\in
  \Omega_k$ reads
\[
[V]_x=\min_{u\in U}\{\beta I[V](x+hf(x+u))+1-\beta\}\,,\;\text{where }\beta=e^{-h}\,,\quad
x \in \Omega_k,
\]
see \cite{FalconeFerretti:2014}, and the Hamiltonian takes the simplified form
\[
H_{\cI}(x,V,u,\cI)=\underset{u\in U_{\cI}}{\min}\{\tilde c_{\cI}u_1+\tilde d_{\cI}u_2+\tilde e_{\cI}\}.
\]
}
$\square$\end{remark}



\noindent To evaluate the right hand side in \eqref{nh1} numerically, we compute \eqref{nh2} for $\mathcal{I}
\in \mathcal W$ by applying numerical optimization methods and then determine the minimizer on the
macro-cell by comparison over the sets $U_{\cI}$.

To determine the minimizer of the right hand side in \eqref{H} we define
\[
F(u)= \frac{\tilde a_{\cI}}{2}u_1^2+\frac{\tilde b_{\cI}}{2}u_2^2+\tilde
c_{\cI}u_1+\tilde d_{\cI}u_2+\tilde e_{\cI}.
\]
Then the minimization problem is given by
\begin{align}\label{minimizing_problem}
  \min_{u\in U_{\cI}} F(u),\quad \mathcal I \in W.
\end{align}

In particular, the following types of cost functionals are of interest for scalars $a,b,c,d,e,r,l,s \in \R$,
where $a,c \ge 0$. Thereby we include also a bilinear term $bu_1u_2$ to allow also bilinear interpolants over quadrangular cells. For the triangular interpolation used throughout this paper, we take $b=0$.

    \begin{enumerate}
      \item Infinite horizon problem with quadratic control cost:
        \begin{align}\label{cost_2}
        F_2(u_1,u_2)=\frac{a}{2}u_1^2+ b u_1 u_2 + \frac{c}{2} u_2^2 + d u_1 + e
      u_2+r.
    \end{align}
    The functional is smooth and convex if $ac-b^2>0$.
    \item Minimum time problem:
      \begin{align}\label{mt}
      F_{\text{MT}}(u_1,u_2)= b u_1 u_2 + d u_1 + e u_2+r.
    \end{align}
    The functional is non-convex if $b\neq 0$ but smooth.
    \item Minimum time/Infinite horizon problem with $\ell_1$-control cost:
      \begin{align}\label{functional_MT_IH}
        F_{\text{MT/IH}}(u_1,u_2)=a|u_1| + b u_1 u_2 + c |u_2| + d u_1 + e u_2+r.
      \end{align}
      The functional is non-convex if $b\neq 0$ and non-smooth.
    \item Infinite horizon problem with $\ell_2$- and $\ell_1$-control cost:
      \begin{align}\label{cost_ih}
      F_{\text{IH}}(u_1,u_2)= \frac{a}{2}u_1^2 + l|u_1| + b u_1 u_2 + \frac{c}{2}u_2^2+ s |u_2| + d
      u_1 + e u_2+r.
    \end{align}
    The functional is convex if $ac-b^2>0$ and non-smooth.
  \end{enumerate}
  \begin{remark}
    { \em Since $U$ is compact all four functionals allow a minimum. If $ac -b^2>0$ it is unique for
    \eqref{cost_2} and \eqref{mt}. For a discussion of optimal control problems of dynamical
systems for functionals \eqref{functional_MT_IH} and \eqref{cost_ih} we refer to
\cite{AltSchneider:2013,MaurerVossen:2006}.}
$\square$\end{remark}


\subsection{Special case: eikonal dynamics}\label{sec:eikonal}
We consider the relation between the subdivision of the state and control space in the context
of minimum time problems with eikonal dynamics. To set up the control problem we introduce a closed target
$\mathcal T \subset \R^n$ , with  $\operatorname{int}(\mathcal T)\neq \emptyset$, and 
smooth boundary.  We consider the minimum time function
\begin{align}
t(x,\tilde u)=\left\{
\begin{array}{ll}
  \inf \left\{ t\in \R_+ : y(t,\tilde u) \in \mathcal T \right\} &
  \text{if } y(t,u(t)) \in \mathcal T\text{ for some } t, \\
  + \infty & \text{otherwise},
\end{array}
\right.
\end{align}
where $\tilde u \in \mathcal U$ and $ y(\cdot ,u( \cdot))$ denotes the solution of
\eqref{eq:dynconaut} depending on the control. Furthermore we assume a small-time local
controllability assumption as formulated in \cite[p. 216]{FalconeFerretti:2014}.
The minimum time problem is defined as 
\[
T(x)=\inf_{\tilde u \in \mathcal U} t(x,\tilde u)
\]
and the minimum time function can be characterized as the viscosity
solution of 
\begin{equation}
\begin{aligned}
  \sup_{u\in U} \left\{ -f(x,u)^T \nabla v(x)\right\}&=1 & &\text{in }\mathcal R\setminus \mathcal T,\\
  v(x)&=0&&\text{on }\partial \mathcal T,
\end{aligned}
\end{equation}
where $\mathcal R$ are all points in the state space for which the time of arrival is finite.

We consider the dynamics 
\begin{equation}\label{eik}
f(x,u)=\left(\begin{array}{c}u_1\\u_2\end{array}\right)\, ,
\quad
U = \{u \in \R^2:\norm{u}_2 \le 1\}, \quad x \in \R^2,
\end{equation}
which leads to an equation of eikonal type.
For this problem we have a direct relation between the direction of the control $u$
and the identification of the arrival point in the state space; for instance, all the arrival points in the first
quadrant will correspond to controls belonging to the set $U_{\{1,2\}}=\{u\in U \; | \;0\leq u_1\,, 0\leq u_2\,\}$.
It is possible to express the interpolation operator \eqref{SLih} for given
node $x$ piecewise as
\begin{equation}
  I[V](x+hf(x,u))=I_{\cI}[V](x+hf(x,u))
\end{equation}
with $x+hf(x,u)$ in quadrant $U_{\cI}$.
For every quadrant we obtain control-dependent formulas of the form
\begin{equation}
  I_{\cI}[V](x+hf(x,u))=c_{\cI} u_1+d_{\cI} u_2+e_{\cI}\,
\end{equation}
with coefficients $c_{\cI},d_{\cI}, e_{\cI}\in \R$.

  \section{Algorithms for solving minimization problems of the form \eqref{minimizing_problem}}\label{sec:alg}
We consider the following approaches:
    \begin{enumerate}
          \item[a)] Minimization by comparison over a finite subset $U_{\text{finite}}\subset U_{\cI}$.
      \item[b)] First order primal-dual method: Chambolle-Pock algorithm, see \cite{ChambollePock:2011}.
      \item[c)] Second-order method: Two different types of semismooth Newton methods depending on the
        smoothness of the cost functional.
      \item[d)] If the functional is of type \eqref{mt}, the controls are located on the
        sphere or in the origin and can be found by a classical Newton method with a suitable chosen
        initialization over the parameterized sphere.

    \end{enumerate}
    
The minimization by comparison approach a) is broadly used in the literature. It consists in
choosing a finite subset $U_{\text{finite}} \subset U$ where the cost function is evaluated and the
minimum is selected among the corresponding values. Such a procedure induces a different
optimization paradigm, in the sense that the continuous nature of the control space is replaced by a
discrete approximation. If a parameter-dependent discretization of the control set is considered,
where errors with respect to the continuous set can be estimated, the discretization has to be fine
enough such that the error is negligible compared to the errors introduced in the different
discretization steps of the overall scheme. Furthermore, accurate discretization of the control
space are far from trivial even in simple cases, as in the three-dimensional eikonal dynamics with
$U=\{u\in\R^3 : \|u\|_2\leq 1\}$, where a spherical coordinate-based discretization of the control
introduces a concentration of points around the poles.  In contrast to this approach, we propose
several numerical algorithms for the treatment of problems of the form \eqref{minimizing_problem}, which preserve the continuous nature of the optimization problem at a similar computational cost.

    \subsection{The smooth case: Cost functionals of type \eqref{cost_2} and
    \eqref{mt}}\label{sec:smooth} In this section we consider control constraints of the type
    $U=\{ u \in \R^m : \norm{u}_2 \le 1\}$.
    To solve the optimization problem for smooth cost functionals of type \eqref{cost_2} and
    \eqref{mt} we present a first-order primal-dual
    algorithm and a semismooth Newton method.

    \subsubsection{Primal-dual algorithm}
    Here we assume that problem \eqref{minimizing_problem} can be reformulated as
    \begin{align}\label{mp_modified}
      \min_{u\in \R^m} F(u)+I_K(u),
    \end{align}
    where $F$ is smooth, and convex and $I_K(u)$ is the indicator function of
    $K=U_{\cI}$.
    This fits in the setting presented in \cite{ChambollePock:2011}, where a primal-dual algorithm (also known as Chambolle-Pock algorithm) is
    formulated, when using the same $F$, and by choosing the
    mappings denoted by $K$ and $G$ in the reference by $\operatorname{id}$
    and $I_K$. We recall the algorithm in the following.

    \vspace{0.1cm}
\begin{algorithm}[H]\label{cp1}\caption{Chambolle-Pock algorithm}
  \SetAlgoLined
  \vskip 1mm
  \KwData{Choose $n=0$, $\tau>0$, $\sigma>0$, $\eta>0$, $\vartheta \in [0,1]$, $(u_0,y_0)\in \R^m
  \times \R^d$, $\bar{u}_0=u_0$.}
  \Repeat{$\norm{x_n - x_{n-1}}_{\R^d\times \R^m \times \R^m }
  < \eta $}
  {  {\vskip 1mm  
  \textsl{Compute $ x_{n+1}=(y_{n+1},u_{n+1},\bar u_{n+1})$ by}
  \begin{equation*}   
    \left\{
    \begin{aligned}
      y_{n+1} &= (I+\sigma \D F^*)^{-1} (y_n+ \sigma \bar{u}_n),\\
      u_{n+1} &= (I+\tau \D G)^{-1}(u_n- \tau y_{n+1}),\\
      \bar{u}_{n+1}&=u_{n+1}+\vartheta (u_{n+1}-u_n).\\
    \end{aligned}
    \right.
  \end{equation*}\vskip -4mm}
  {
  \vskip 1mm
  \textsl{Set $n=n+1$}.
  }
  \vskip 1mm
  }
\end{algorithm}
    \vspace{0.1cm}

    From identities from convex analysis, see \cite{Rockafellar:1997}, we have
    \begin{align}
      & (I+\tau \D G)^{-1}(y)=\underset{\operatorname{u \in \R^m}}{\operatorname{argmin}} \left\{
      \frac{\norm{u-y}_2^2}{\tau}+I_K(u) \right\}=P_K(y) ,\\\label{identity3}
      & u-(I+\tau \D F)^{-1}(u)=\tau \left(I+\frac{1}{\tau}\D F^*\right)^{-1}\left(\frac{u}{\tau}\right),
    \end{align}
    where $P_K$ denotes the projection on $K$.
    From \eqref{identity3} we have with $\sigma=\frac{1}{\tau}$
        \begin{align*}
          (I+\sigma \partial F^*)^{-1} (\sigma u)&=\sigma u - \sigma \left(I + \frac{1}{\sigma} \D
          F\right)^{-1}(u), 
        \end{align*}
        and hence
        \begin{align*}
          (I+\sigma \partial F^*)^{-1} (u)&= u - \sigma \left(I + \frac{1}{\sigma} \D
          F\right)^{-1}\left(\frac{1}{\sigma}u\right).
        \end{align*}

        In the case $\cI=\{1,2\}$ the projection $P_K$ is given by 
\begin{align}\label{projection}
P_K\colon \R^m \rightarrow \R^m,\quad P_K(p)=\frac{\max(0,p)}{\max(1,\norm{\max(0,p)}_2)} .
\end{align}
%
This leads to Algorithm \ref{cp2}.

    \vspace{0.1cm}
\begin{algorithm}[H]\label{cp2}\caption{Chambolle-Pock algorithm (variant)}
  \SetAlgoLined
  \vskip 1mm
  \KwData{Choose $n=0$, $\tau>0$, $\sigma>0$, $\eta>0$, $\vartheta \in [0,1]$, $(u_0,y_0)\in \R^m
  \times \R^d$, $\bar{u}_0=u_0$. 
  }
  \Repeat{ $\norm{x_n - x_{n-1}}_{\R^d\times \R^m \times \R^m }
  < \eta $}
  {  {\vskip 1mm  
  \textsl{Compute $x_{n+1}=(y_{n+1},u_{n+1},\bar u_{n+1})$ by}
          \begin{equation*}
            \left\{
            \begin{aligned}
              y_{n+1} &= y_n+ \sigma\bar{u}_n - \sigma \left(I + \frac{1}{\sigma} \D
              F\right)^{-1}\left(\frac{y_n}{\sigma}+ \bar{u}_n\right),\\
              u_{n+1} &= P_K(u_n - \tau y_{n+1}),\\
              \bar{u}_{n+1}&=u_{n+1}+\vartheta (u_{n+1}-u_n).\\
            \end{aligned}
            \right.
          \end{equation*}
\vskip -4mm}
  {\vskip 1mm
  \textsl{Set $n=n+1$}.
  }}
\end{algorithm}
    
\vspace{0.1cm}
      If $(I+\frac{1}{\sigma} \D F)^{-1}$ is linear, the first step can be reformulated as
      \begin{equation*}
        \begin{aligned}
          y_{n+1} &= \left(I+\frac{1}{\sigma} \D F\right)^{-1} \D F\left(\frac{y_n}{\sigma}+ \bar{u}_n\right).
        \end{aligned}
      \end{equation*}
    \begin{remark}
{\em Since the cost functional has only quadratic and linear terms it is more convenient to use a linear
interpolant rather than a bilinear one in the semi-Lagrangian scheme.This does not produce any bilinear terms. 
If $F$ has a bilinear term it is very challenging to compute $(I +
\frac{1}{\sigma} \D F)^{-1}$ in higher dimensions, since then $ \D F$ depends nonlinearly on $u$.}
    $\square$\end{remark}

    \subsection{Semismooth Newton method}
The Chambolle-Pock algorithm is a first order algorithm which requires convexity of the functional.
We refer to \cite{Valkonen:2014} for an extension to a class of non-convex problems. 
As a second algorithm that we propose we turn to the semismooth Newton method which does not rely on global
convexity. We recall some main aspects of semismooth functions, cf. \cite{Ulbrich:2011}.
\begin{definition}[Semismoothness]
  Let $V\subset \R^{\mu}$ be nonempty and open, $\mu \in \N$. Then function $f\colon V \rightarrow
  \R^{\mu}$ is semismooth
at $x\in V$ if it is Lipschitz continuous near $x$ and if the following limit exists for all $s\in
\R^{\nu}$, $\nu \in \N$:
\[
\lim_{M\in \partial f(x+\tau d)_{d \rightarrow s, \tau \rightarrow 0^+} } Md.
\]
\end{definition}
Furthermore, there holds the following chain rule.
\begin{lemma}[Chain rule]
  Let $V\subset \R^{\mu}$ and $W\subset \R^{\nu}$ be nonempty open sets, $g\colon V\rightarrow W$ be
  semismooth at $x\in V$, and $h\colon W \rightarrow \R^{\eta}$, $\eta \in \N$, be semismooth at $g(x)$ with $g(V)\subset
W$. Then the composite map $f:=h\circ g \colon V \rightarrow \R^{\eta}$ is semismooth at $x$. Moreover,
\[
f'(x,\cdot)=h'(g(x),g'(x,\cdot)).
\]
\end{lemma}

To set up a semismooth Newton algorithm for \eqref{minimizing_problem} we proceed as follows. For simplicity, we focus on
the case $\cI=\{1,2\}$. The first-order optimality condition
for cost functionals of type \eqref{cost_2} and \eqref{mt} can be formulated as 
\begin{align}\label{opt_cond}
  u=P_K(u-\vartheta \nabla F(u))\quad \forall \vartheta >0
\end{align}
for $u\in \R^m$, $F=F_2$ or $F=F_{\text{MT}}$
with projection $P_K$ as in \eqref{projection}.
We introduce $p=y-\vartheta \nabla F(y)$ and choose $\vartheta$ such that $\vartheta \nabla F(u)
$ is of the same scale as $y$ and rewrite \eqref{opt_cond} as
\begin{align*}
  u-P_K(p)&=0,\\
  u-\vartheta \nabla F(u)-p&=0.
\end{align*}
Setting $\beta=\max(1,\norm{\max(0,p)}_2)$ we define
\[
E(u,p,\beta)=
\begin{pmatrix}
  \beta u - \max(0,p) \\
  u-\vartheta \nabla F(u)-p\\
  \beta-\max(1,\norm{ \max(0,p)}_2) \\
\end{pmatrix}.
\]
Now the optimality condition can be formulated as 
\[
E(z)=0
\]
with $z=(u,p,\beta) $.
Then we set up a semismooth Newton method as 
presented in Algorithm \ref{ssn1}.

    \vspace{0.1cm}
\begin{algorithm}[H]\label{ssn1}\caption{Semismooth Newton algorithm}
  \SetAlgoLined
  \vskip 1mm
  \KwData{ Choose $n=0$, initialization $z_0 \in U_{\cI} \times \R^m \times \R$, $\eta>0$.
  }
  \Repeat{$\norm{z_n - z_{n-1}}_{\R^m\times \R^m \times \R  }< \eta $}
  {  {\vskip 1mm  
  \textsl{Solve the Newton equation for $\delta z\in U_{\mathcal
  I}\times \R^m \times \R$ } given by 
  \begin{align}\label{SSN_Iteration_matrix}
    \begin{pmatrix}
      \beta_n I & D_{\chi_{p_n\ge 0}}& y_n \\[1ex]
      I-\vartheta\nabla^2 F(u_n) & - I & 0 \\[1ex]
      0 & \chi_{m_n\ge 1}\frac{-p_n^T D_{\chi_{p_n\ge 0}}}{m_n} & 1 \\
    \end{pmatrix}      \delta z
    = -E(z_n)
  \end{align}
 (with matrix $D_{\chi_{p_n\ge 0}}=\operatorname{diag}(\chi_{p_n\ge 0})$ and $ m_n=
 \norm{\max(0,p_n)}_2$).
 {\vskip 1mm
\textsl{Update}
    \begin{align}
    z_{n+1}=z_n + \delta z.
  \end{align}}
  \vskip -4mm}
  {\vskip 1mm
  \textsl{Set $n=n+1$}.
  }}
\end{algorithm}

\vspace{0.1cm}

We address solvability of \eqref{SSN_Iteration_matrix}.
\begin{lemma}\label{lem:regularity}
Set $M=  I-\vartheta\nabla^2 F(u_n)$ and $q^+=\max(0,q)$ for $q\in \R^m$. The Newton matrix
is regular if one of the following two conditions is satisfied:
\begin{enumerate}
  \item $m_n< 1$: $\beta_n\neq 0$, $ \left(\frac{1}{\beta_n} M D_{\chi_{p_n}\ge 0}   + I
    \right)$ regular.
\item $m_n\ge 1$: $\beta_n\neq 0$, 
  \begin{align}\label{regularity_cond}
  I -\frac{1}{\beta_n} M\quad  \text{regular}
\end{align}
  and\quad 
  $
  \left\{-\frac{1}{m_n} (p^+)^T (I -\frac{1}{\beta_n} M^+)^{-1} \left( - M^+  \frac{u^+_n}{\beta_n}
  +z_4^+) \right)^+ +1 \right\}\neq 0
  $
  , where the components of $p=(p^+,p^-)$ are ordered with respect to active and inactive sets and
  $M^+$ denotes the submatrix corresponding to $p^+$. 
  \end{enumerate}
\end{lemma}

\begin{remark}
  { \em For $\beta_n \approx 1$ condition \eqref{regularity_cond} is closely related to the regularity of
  the Hessian of $F$. Moreover, if we set $K=\sup_{u\in U_{\{1,2\}}} \norm{\nabla^2F(u)}$ and assume that $\beta_n$
  are bounded away from $1$ (i.e. $\beta_n\ge \kappa >1$) then the choice
  $\vartheta<\frac{\kappa-1}{K}$ implies \eqref{regularity_cond}.
}
$\square$\end{remark}

{\em Proof of Lemma \ref{lem:regularity}}.
The two cases are considered separately.
For $m_n < 1$, the Newton-Matrix is given by
    \begin{align*}
DE(z_n)=\begin{pmatrix}
  \beta_n I & D_{\chi_{p_n\ge 0}} & u_n \\
  M & - I & 0 \\
  0 & 0 & 1 \\
\end{pmatrix}.
    \end{align*}
The regularity of the matrix follows from the regularity of $\begin{pmatrix}
  \beta_n I & D_{\chi_{p_n\ge 0}} \\
  M & - I 
\end{pmatrix}$ and hence of
$
\left(\frac{1}{\beta_n} M D_{\chi_{p_n}\ge 0}   + I \right).
$

For $m_n \ge 1$ the Newton-matrix is given by
  \begin{align*}
    DE(z_n)=\begin{pmatrix}
      \beta_n I & D_{\chi_{p_n\ge 0}}& u_n \\
 M & - I & 0 \\
  0 & \frac{-(p^+)^T}{m_n} & 1 
\end{pmatrix}
\end{align*}
and leads to the Newton equations 
\begin{align}
\beta_n \delta u -\delta p^+  + \delta \beta u_n &=z_1,\\
M \delta u - \delta p &= z_2 , \\\label{Eq3}
-\frac{1}{m_n} (p^+)^T \delta p^+ + \delta \beta &= z_3
\end{align}
for suitable $z_1,z_2 \in \R^m$, $z_3 \in \R$ independent of $(\delta u, \delta p, \delta  \beta)$.
From the first and second equations we obtain
\begin{align*}
  -\frac{1}{\beta_n} M (\delta p^+ - \delta \beta u_n) + \delta p &=z_4
\end{align*}
for $z_4 \in \R^m$.
Without loss of generality we assume that the components of $\delta p=(\delta p^+, \delta p^-)$  are ordered with
respect to active
and inactive sets. Then we have
\begin{align}
  (I -\frac{1}{\beta_n} M^+) \delta p^+ + \delta \beta M^+  \frac{u^+_n}{\beta_n} &=z_4 
\end{align}
and therefore
\begin{align}
  \delta p^+ =  \delta \beta (I -\frac{1}{\beta_n} M^+)^{-1} \left( - M^+  \frac{u^+_n}{\beta_n}
  +z_4^+\right). 
\end{align}

With \eqref{Eq3} we obtain
\begin{align}
  \left\{-\frac{1}{m_n} (p^+)^T (I -\frac{1}{\beta_n} M^+)^{-1} \left( - M^+  \frac{u^+_n}{\beta_n}
  +z_4^+) \right)^+ +1 \right\}\delta  \beta =\bar \beta.
\end{align}
\begin{flushright}
$\square$
\end{flushright}

\begin{theorem}\label{thm1}
    Let $\bar{u}$ be a strict local minimizer of \eqref{minimizing_problem}. Under the assumptions of Lemma~\ref{lem:regularity} the semismooth Newton method converges locally
    superlinearly or terminates after a finite number of steps at $\bar{u}$.
\end{theorem}

\proof
 The operator $E$ is semismooth and is as a composition of Lipschitz continuous functions again
 Lipschitz continuous. Furthermore from Lemma \ref{lem:regularity} we obtain the boundedness of the
 inverse derivative in a neighbourhood of $\bar{u}$.
 
 Consequently, the assertion follows from \cite[p. 29]{Ulbrich:2011} and \cite[p. 220, Thm.
 8.3]{ItoKunisch:2008}.
\endproof


\subsection{Approach for cost functionals of type \eqref{mt}}

To treat the problem \eqref{minimizing_problem} for cost functionals of type \eqref{mt} we present
an alternative approach. 
We show that all possible minimizers are
located on the sphere or in the origin; cf. also the discussion in \cite{CristianiFalcone:2007} where a linear functional is considered. To minimize over the sphere we parameterize the sphere by polar coordinates and consider the restriction
of the functional on the sphere.
There the minimizer can be found
by applying a classical Newton method.

\begin{lemma}
  For control problems \eqref{minimizing_problem} with cost functional \eqref{mt} for $r=0$ all minimizers are located either on the sphere or in the origin.
\end{lemma}
\proof
  For simplicity we consider the case \eqref{minimizing_problem} for $d=2$ and $\mathcal I=\{1,2\}$.
  Let $F(u)=du_1+bu_1u_2+eu_2$, $u\in \R^2$, $d,b,e \in \R$. Then we have
\[
\nabla F(u)= \begin{pmatrix}
  d+bu_2\\
  e+bu_1
\end{pmatrix},
\quad
\nabla^2 F(u)=\begin{pmatrix} 0 & b \\ b & 0 \end{pmatrix}.
\]
This implies, we have for $b\neq 0$ a saddle point with eigenvalues $\pm 1$ in
\[
u_1=-\frac{e}{b},\quad u_2=-\frac{d}{b}.
\]
For $b=0$ the minimum is reached in $u=(0,0)^T$.
Consequently, all minimizers are on the boundary of $\partial U_{\{1,2\}}$. However, since the
restriction of the bilinear functional to the axes is linear the assertion follows immediately.
\endproof

In two dimensions the functional $F$ restricted to the sphere is given by
\[
\tilde F(\varphi)=a\cos (\varphi) + \frac{b}{2} \sin (2\varphi) + c \sin(\varphi),\quad 0\le \varphi
\le
\frac{\pi}{2}.
\]



%

%
%
\subsection{The non-smooth case: Cost functionals of type \eqref{functional_MT_IH} and \eqref{cost_ih}}
In this section we consider two types of control constraints, box constraints as well as Euclidean
constraints. Finally,  we also present a  splitting approach.

%

\subsubsection{Semismooth Newton in case of Euclidean norm constraints}\label{sec:ss}

In this section we restrict the consideration to equations of eikonal type and present the approach for a two
dimensional problem. Therefore the control set is given by $U=\{u\in\R^ 2: \|u\|_2\leq 1\}$.
%
Let us consider the equivalent problem formulation
\begin{align}\label{nonsmooth_problem}
  \min_{u \in \R^2} E(u)+h(u),
\end{align}
where $E\colon \R^2 \rightarrow \R$ is smooth and 
\begin{equation}
h(u)= \left\{
\begin{array}{ll}
  \alpha( |u_1|+|u_2|),  &\text{if }u \in U_{\cI},\\
  \infty, & \text{if }u \not \in U_{\cI}.
\end{array}
\right.
\end{equation}
The optimality condition for the nonsmooth problem \eqref{nonsmooth_problem} is given by
\[
0 \in \nabla E(u)+\partial h(u)\,,
\]
where $\partial h$ denotes the subdifferential of the convex function $h$. Setting $q=-\nabla E(u)$ the optimality condition can be equivalently written as
\begin{equation}
  \left\{
\begin{aligned}
  &q=-\nabla E(u),\\
  &q \in \partial h(u).
\end{aligned}
\right.
\end{equation}
With the convex conjugate we obtain equivalently
\begin{equation}
  \left\{
\begin{aligned}
  &q=-\nabla E(u),\\
  &u \in \partial h^*(q)\,,
\end{aligned}
\right.
\end{equation}
with
\[
h^*(q)=\sup_{v \in U_{\cI}}(v\cdot q -h(v)). 
\]
Note that on $U_{\cI}$
\[
v\cdot q-h(v)=
  (q_1-\alpha_1)v_1+(q_2-\alpha_2)v_2
\]
with
\begin{align}\label{alpha1_alpha2}
\alpha_1=\operatorname{sgn}(L_1) \alpha,\quad \alpha_2=\operatorname{sgn}(L_2) \alpha
\end{align}
and
\begin{align}\label{g_1_g_2}
L_1=\left\{
\begin{array}{ll}
  1, &\text{if } \cI=\{1,2\} \text{ or }\cI=\{1\},\\
  -1, &\text{if } \cI=\{2\} \text{ or }\cI=\emptyset,\\
\end{array}
\right.,
\quad L_2=\left\{
\begin{array}{ll}
  1, &\text{if } \cI=\{1,2\} \text{ or }I=\{2\},\\
  -1, &\text{if } I=\{1\}\text{ or }\cI=\emptyset .
\end{array}
\right.
\end{align}
For fixed $q$ and $v\in U_{\cI}$ we define
\[
l(v)=(q_1-\alpha_1)v_1+(q_2-\alpha_2)v_2.
\]
The $\sup$ of this linear function is attained on $\partial U_{\cI}$. We distinguish four cases
\begin{equation*}
\begin{aligned}
  \sgn(q_1-\alpha_1)\neq \sgn(L_1)\;\land\;\sgn(q_2-\alpha_2)\neq \sgn(L_2):\quad & h_{L_1,L_2}^*(q)=0,\\
\sgn(q_1-\alpha_1)=\sgn(L_1)\;\land\; \sgn(q_2-\alpha_2)\neq \sgn(L_2):\quad  & h_{L_1,L_2}^*(q)=q_1 -
\alpha,\\
\sgn(q_1-\alpha_1)\neq \sgn(L_1)\;\land\; \sgn(q_2-\alpha_2)=\sgn(L_2):\quad & h_{L_1,L_2}^*(q)=q_2 - \alpha,\\
\sgn(q_1-\alpha_1)=\sgn(L_1)\;\land\; \sgn(q_2-\alpha_2)=\sgn(L_2):\quad & \\
h_{L_1,L_2}^*(q)= \sup_{\varphi \in
\Lambda} (q_1 - \alpha_1) &\cos \varphi  +(q_2 - \alpha_2) \sin \varphi\,,
\end{aligned}
\end{equation*}
with
\begin{align}
  \Lambda=\left\{
  \begin{array}{ll}
    0\le \varphi \le \frac{\pi}{2}& \text{if }\cI=\{1,2\},\\
    \frac{\pi}{2}< \varphi \le \pi& \text{if }\cI=\{2\},\\
    \pi < \varphi \le \frac{3\pi}{2}& \text{if }\cI=\emptyset,\\
    \frac{3\pi}{2}< \varphi \le 2\pi & \text{if }\cI=\{1\}.
  \end{array}
  \right.
\end{align}

In the following for simplicity we only consider the problem for $\cI ={\{1,2\}}$.
From the first order condition we obtain the following relation between $q$ and $\varphi$
\begin{align}\label{relation_q_phi}
\tan \varphi= \frac{q_2 - \alpha}{q_1 - \alpha} >0.
\end{align}

Summarizing we have
\begin{equation}
  h^*(q)=
  \left\{
  \begin{aligned}
    &(q_1 - \alpha) \cos \varphi + (q_2 - \alpha) \sin \varphi,& & \text{if }q_1,q_2\ge \alpha, \\
& q_2 - \alpha, & &\text{if } q_1<\alpha, \quad q_2>\alpha, \\
&0,& & \text{if }q_1,q_2<\alpha, \\
&q_1 - \alpha ,& & \text{if }q_1>\alpha ,q_2<\alpha\,,
  \end{aligned}
  \right.
\end{equation}
with $\varphi$ given by \eqref{relation_q_phi} and for the generalized derivative

\begin{equation}
 \partial h^*(q)=
  \left\{
  \begin{aligned}
    & (w_1(q),w_2(q))^T,& & \text{if }q_1,\; q_2>\alpha, \\
    & (0, 1 )^T, && \text{if } q_1<\alpha,\; q_2>\alpha, \\
&(0, 0)^T,& & \text{if }q_1<\alpha,\; q_2<\alpha, \\
& (1 , 0)^T,& & \text{if }q_1>\alpha,\; q_2<\alpha\,
  \end{aligned}
  \right.
\end{equation}
where
\[(w_1(q),w_2(q))=\nabla((q_1 - \alpha) \cos \varphi + (q_2 - \alpha)\sin \varphi)\,.\]


To set up a semismooth Newton method we introduce a piecewise affine, continuous regularization of
the gradient in a $\varepsilon$-tube around the discontinuity
and define

\begin{equation}
  (\partial h^*)_{\varepsilon}(q)=
  \left\{
  \begin{aligned}
    & (w_1(q),w_2(q))^T,& & \text{if }q_1,\; q_2>\alpha+\varepsilon, \\
    & (0, 1 )^T, && \text{if } q_1<\alpha,\; q_2>\alpha+\varepsilon, \\
        & \left(0, \frac{q_2-\alpha}{\varepsilon} \right)^T, && \text{if } q_1<\alpha,\; \alpha<q_2<\alpha+\varepsilon, \\
&(0, 0)^T,& & \text{if }q_1<\alpha,\; q_2<\alpha, \\
        & \left(\frac{q_1-\alpha}{\varepsilon},0 \right)^T, && \text{if } \alpha<q_1<\alpha+\varepsilon,\; q_2<\alpha, \\
& (1 , 0)^T,& & \text{if }q_1>\alpha+\varepsilon,\; q_2<\alpha,\\
& \frac{r}{\varepsilon}(w_1(q_p) , w_2(q_p))^T,& & \text{if }\alpha<q_1,\; \alpha<q_2,\:\text{and }
\|q-(\alpha,\alpha)^T\|_2\leq\varepsilon,
  \end{aligned}
  \right.
\end{equation}
where 
\[r=\|q-(\alpha,\alpha)\|_2\,,\quad\theta=\tan^{-1}\left(\frac{q_2-\alpha}{q_1-\alpha}\right)\,,\;\text{and
}\; q_p=\varepsilon(\cos(\theta),\sin(\theta))+(\alpha,\alpha).\]
%

The semismooth Newton method is presented in Algorithm \ref{ssn2}.

\vspace{0.1cm}
\begin{algorithm}[H]\label{ssn2}\caption{Semismooth Newton algorithm}
  \SetAlgoLined
  \vskip 1mm
  \KwData{ Choose $n=0$, initialization $z_0=(q_0,u_0) \in \R^2 \times U_{\cI}$, $\eta>0$.
  }
  \Repeat{ $\norm{z_n - z_{n-1}}_{\R^2\times \R^2}< \eta $}
  {  {\vskip 1mm  
  \textsl{Solve the Newton equation for $\delta z\in \R^2 \times U_{\mathcal{I}}$ given by}
  \begin{align}\label{SSN_2}
     \begin{pmatrix}
  I & \nabla^2 E(u_n) \\ -D(\partial h^*)_{\varepsilon} (q_n)& I 
\end{pmatrix}
  \delta z
    = \begin{pmatrix}
q_n + \nabla E(u_n) \\
- (\partial h^*)_{\varepsilon}(q_n) + u_n
\end{pmatrix},
  \end{align}
  {\em  where $ D(\cdot)$ denotes the Newton derivative.}
 {\vskip 1mm
\textsl{Update}
    \begin{align}
    z_{n+1}=z_n + \delta z.
  \end{align}}
  \vskip -4mm}
  {\vskip 1mm
  \textsl{Set $n=n+1$}.
  }}
\end{algorithm}
    \vspace{0.1cm}

    We can conclude the following theorem.
\begin{theorem}\label{thm:SSN_2} Let $u \in \R^2$ be a strict local minimum of \eqref{nonsmooth_problem}.
If
\begin{align}\label{cond_invert}
I-\nabla^2 E(u_n) D(\partial h^*)_{\varepsilon}(q_n)
\end{align}
is regular for all $k$, the semismooth Newton converges locally superlinearly.
\end{theorem}
\proof
We introduce the piecewise affine function 
\[
G\colon \R^2 \times \R^2 \rightarrow \R^2\times \R^2,\quad G(q,u)=(q_n+\nabla E(u),-(\partial h^*
)_{\varepsilon}(q)+u_n).
\]
$G$ is semismooth and together with condition \eqref{cond_invert} we get directly the regularity of the Hessian in
the Newton equation given in \eqref{SSN_2}. From both conditions we derive locally superlinear
convergence as in Theorem \ref{thm1}.
\endproof

\subsubsection{Semismooth Newton in case of box constraints}\label{sec:ssbc}

Now, let the control set be given by $U=\{u\in \R^m : u_a \le u \le u_b\}$ with $u_a\le 0 \le u_b$,
$u_a,u_b \in \R^2$. As in the previous section, we consider a cost functional of the type

\begin{align}
  \min_{u \in \R^2} E(u)+h(u),
\end{align}
where $E$ is smooth and 
\begin{equation}
h(u)= \left\{
\begin{array}{ll}
  \alpha( |u_1|+|u_2|), & \text{if }u \in U_{\cI},\\
  \infty, &\text{if }u \not \in U_{\cI}.
\end{array}
\right.
\end{equation}

\noindent Proceeding as in the previous section we distinguish four cases
\begin{equation*}
\begin{aligned}
  \sgn(q_1-\alpha_1)\neq \sgn(L_1)\;\land\; \sgn(q_2-\alpha_2)\neq \sgn(L_2):\quad & h_{L_1,L_2}^*(q)=0,\\
\sgn(q_1-\alpha_1)=\sgn(L_1)\;\land\; \sgn(q_2-\alpha_2)\neq \sgn(L_2):\quad  &
h_{L_1,L_2}^*(q)=|q_1 - \alpha_1|u_a,\\
\sgn(q_1-\alpha_1)\neq \sgn(L_1)\;\land\; \sgn(q_2-\alpha_2)= \sgn(L_2):\quad &
h_{L_1,L_2}^*(q)=|q_2 - \alpha_2|u_b,\\
\sgn(q_1-\alpha_1)=\sgn(L_1)\;\land\; \sgn(q_2-\alpha_2)= \sgn(L_2):\quad & h_{L_1,L_2}^*(q)=|q_1
-\alpha_1|u_a  \\
&\hspace{1.5cm} + |q_2  - \alpha_2|u_b\,,
\end{aligned}
\end{equation*}
with $\alpha_1,\alpha_2$ as in \eqref{alpha1_alpha2} and $L_1,L_2$ as in \eqref{g_1_g_2}. 
Further, there holds
\begin{equation}
  \nabla h^*(q)=
  \left\{
  \begin{aligned}
    & (\sgn(q_1)\,u_a,\sgn(q_2)\,u_b)^T, && \text{if } q_1>\alpha,\; q_2>\alpha, \\
    & (0,\sgn(q_2)\,u_b)^T, && \text{if } q_1<\alpha,\; q_2>\alpha, \\
&(0, 0)^T,& & \text{if }q_1<\alpha,\; q_2<\alpha, \\
& (\sgn(q_1)\,u_a,0)^T,& & \text{if }q_1>\alpha,\; q_2<\alpha\,,
  \end{aligned}
  \right.
\end{equation}
for $\cI=\{1,2\}$ and accordingly for the other cases. 
To set up a semismooth Newton method we regularize the first
component of the gradient with ramps of width $\varepsilon$ at $|q_1|=\alpha$ and the second component with
ramps at $|q_2|=\alpha$. The overall algorithm has the same form as the semismooth Newton method in
\eqref{SSN_2}.

By a similar consideration as in the previous section Theorem \ref{thm:SSN_2} holds also for the
semismooth Newton method described above.

\subsubsection{A splitting approach}

To treat the cost functionals \eqref{functional_MT_IH} and \eqref{cost_ih} with $\ell_1$-terms we can reformulate the problem without non-differentiable
terms by doubling the number of variables. This is illustrated for the two dimensional case with two
dimensional control $(u,v) \in U_{\cI}$.
For a scalar $z\in \R$ we define $z_+=\max(0,z),\; z_-=\min(0,z)$ (in particular $z_+z_-=0$ and
$|z|=z_+-z_-$).
Then problem \eqref{minimizing_problem} with cost functional \eqref{cost_ih} is given by
\begin{equation*}
\left\{
  \begin{aligned}
    &\min_{u_+,u_-,v_+,v_-} au_+^2+au_-^2+cv_+^2+cv_-^2+bu_+v_+-bu_-v_+-bu_+v_-+bu_-v_- \\
  &\quad +l(u_+-u_-)+s(v_+-v_-) + d(u_++u_-)+e(v_+ + v_-) + r, \quad \text{s.t.}\\
  &g_i+b_i\begin{pmatrix} u_+ - u_- \\ v_+ - v_- \end{pmatrix} \ge 0, \\
  &u_+u_-=0,\quad v_+v_-=0,\\
  &u_+ \ge 0,\quad v_+ \ge 0,\quad u_- \ge 0, \quad v_- \ge 0,\\
  &\norm{u}_2^2+\norm{v}_2^2\le 1, \quad u,v \in U_{\cI}
\end{aligned}
\right.
\end{equation*}
with $g$ chosen as in \eqref{form_dynamics} and $b_i$ as the columns of $B$.
This problem formulation can be simplified, by distinguishing four cases depending on $\cI$;
for example, for $U_{\{1,2\}}$ we obtain the equivalent problem
\begin{equation}\label{problem_splitting}
\left\{
  \begin{aligned}
    &\min_{u_+,v_+} au_+^2+cv_+^2+bu_+v_++(l+d)u_++(s+e)v_+ + r, \quad \text{s.t.}\\
  &g_i+b_i\begin{pmatrix} u_+\\ v_+  \end{pmatrix} \ge 0, \quad i=1,2,\\
  &u_+ \ge 0,\quad v_+ \ge 0,\\
  & u_+^2+v_+^2 \le 1.
\end{aligned}
\right.
\end{equation}

Problem \eqref{problem_splitting} can be solved by using algorithms of Section \ref{sec:smooth}.

%

\section{Numerical examples}\label{sec:examples}
In this section we present a set of numerical tests aiming at studying the performance and accuracy
of the algorithms presented in the previous sections.  We begin by assessing the performance of the proposed numerical optimization routines as a separate building block. 

\subsection{Preliminary tests}
We consider a generic two-dimensional minimization problem of the form
\begin{equation}
\underset{u\in U}{\min} \,\,\frac12 \|u\|_2^2+L\cdot u+\gamma\|u\|_1\
\end{equation}
subject to Euclidean norm or box constraints, i.e.  
$U=\{u \in \R^2 \; |\; \norm{u}_2 \le 1\}$ or 
$U=\{u \in \R^2 \; |\; 0 \leq u_1\leq 1,\,0 \leq u_2\leq 1\}$, respectively.
Results presented in Tables \ref{tab1},  \ref{tab2} and \ref{tab3} show the different performance
scenarios found under different costs and constraints. For every setting we can observe that the
minimization by the comparison algorithm (see Section \ref{sec:alg}) requires very fine discretizations of the control variable (and thus, higher CPU time) to achieve similar error levels as the  optimization-based counterpart.

\begin{table}[htb]
  \begin{tabular}{cccll}
    \hline\\
    Algorithm & Tolerance & Iterations& CPU time  & $\|\cdot\|_2$ error \\
    \cmidrule(r){1-1} \cmidrule(lr){2-2} \cmidrule(lr){3-3} \cmidrule(lr){4-4}\cmidrule(lr){5-5} \\
    Chambolle-Pock     &1E-4 &    9 &   7.3E-5 [s] &  4.31E-5 \\
    Semismooth Newton &1E-4 &    5 &  1.3E-2 [s] &  7.74E-9 \\
    Comparison ($1E4$ evaluations)    & -- & -- &   1.1E-3 [s] &  1.6E-2 \\
    Comparison ($2E3$ evaluations)    & -- & -- &   6.6E-4 [s] &  4.1E-2 \\
    \hline 
  \end{tabular}
  \caption{Performance tests for an $\ell_2$-cost (i.e. $\gamma=0$) with Euclidean norm constraint.}\label{tab1}
\end{table}
\begin{table}[htb]
  \begin{tabular}{cccll}
    \hline\\
   Algorithm & Tolerance & Iterations & CPU time  & $\|\cdot\|_2$ error  \\
    \cmidrule(r){1-1} \cmidrule(lr){2-2} \cmidrule(lr){3-3} \cmidrule(lr){4-4}\cmidrule(lr){5-5} \\
    Semismooth Newton &1E-4 &    101 &  4.53E-3 [s]& 1.51E-3\\
    Comparison ($4E4$ evaluations)    & -- & -- &   5.18E-3 [s] &  5.77E-3  \\
    Comparison ($2E5$ evaluations)    & -- & -- &   2.08E-2 [s] &  2.80E-3\\
    \hline 
  \end{tabular}
  \caption{Combined $\ell_2$- and $\ell_1$-cost ($\gamma=0.1$) with box constraint.}\label{tab2}
\end{table}
\begin{table}[htb]
  \begin{tabular}{cccll}
    \hline\\
    Algorithm & Tolerance & Iterations & CPU time  & $\|\cdot\|_2$ error  \\
    \cmidrule(r){1-1} \cmidrule(lr){2-2} \cmidrule(lr){3-3} \cmidrule(lr){4-4}\cmidrule(lr){5-5} \\
    Semismooth Newton                &1E-4 & 58  &  1.47E-2 [s]& 1.23E-3\\
    Comparison ($2E3$ evaluations)    & --  & -- &   8.56E-4 [s] &  4.18E-2 \\
    Comparison ($1E5$ evaluations)    & --  & -- &   1.20E-2 [s] &  1.47E-2\\
    \hline 
  \end{tabular}
  \caption{Combined $\ell_2$- and $\ell_1$-cost ($\gamma=0.1$) with Euclidean norm constraint.}\label{tab3}
\end{table}

\subsection{Interplay with the semi-Lagrangian scheme}
\begin{figure}[!ht]
  \centering
  \epsfig{file=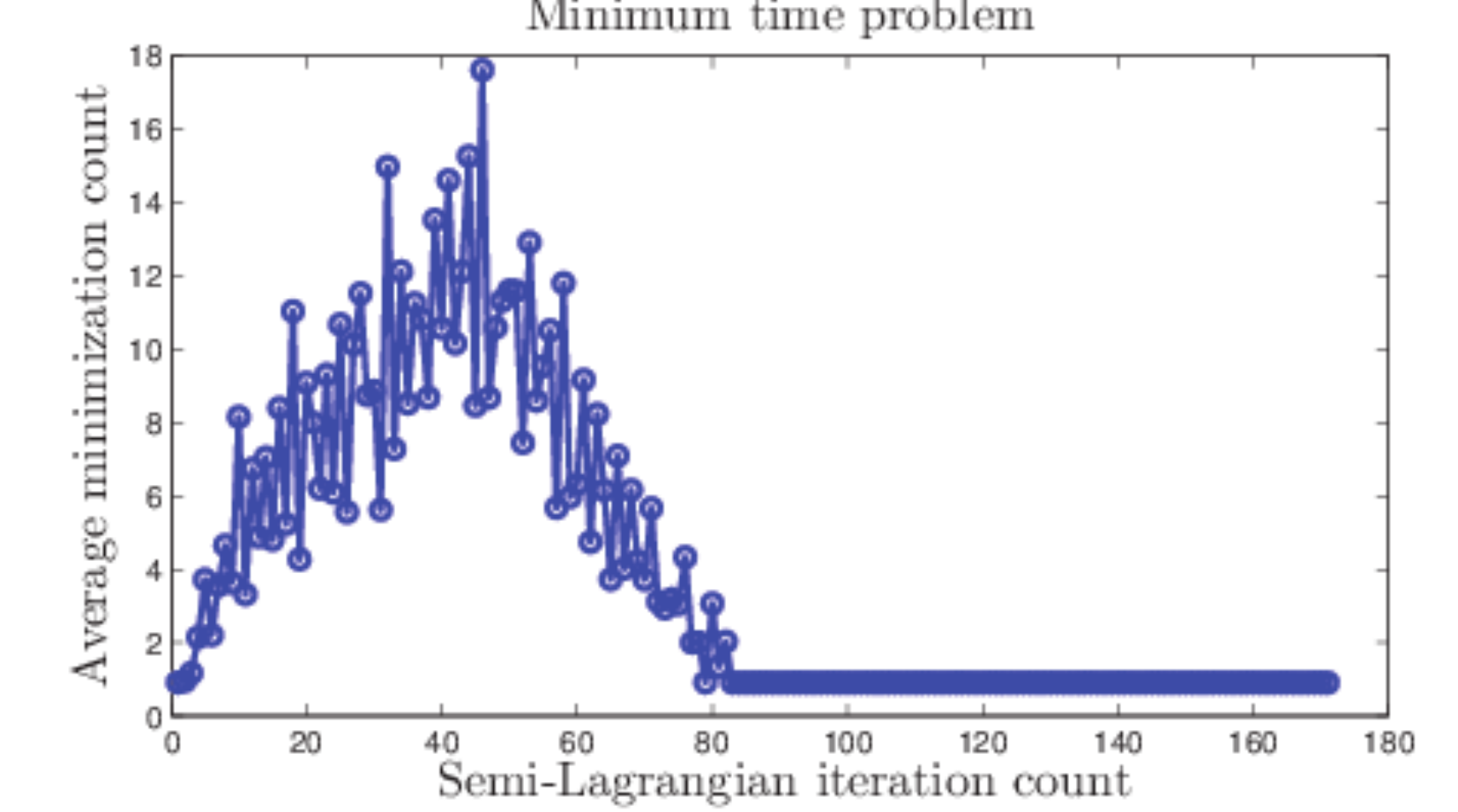,width=.45\linewidth,clip=}
  \epsfig{file=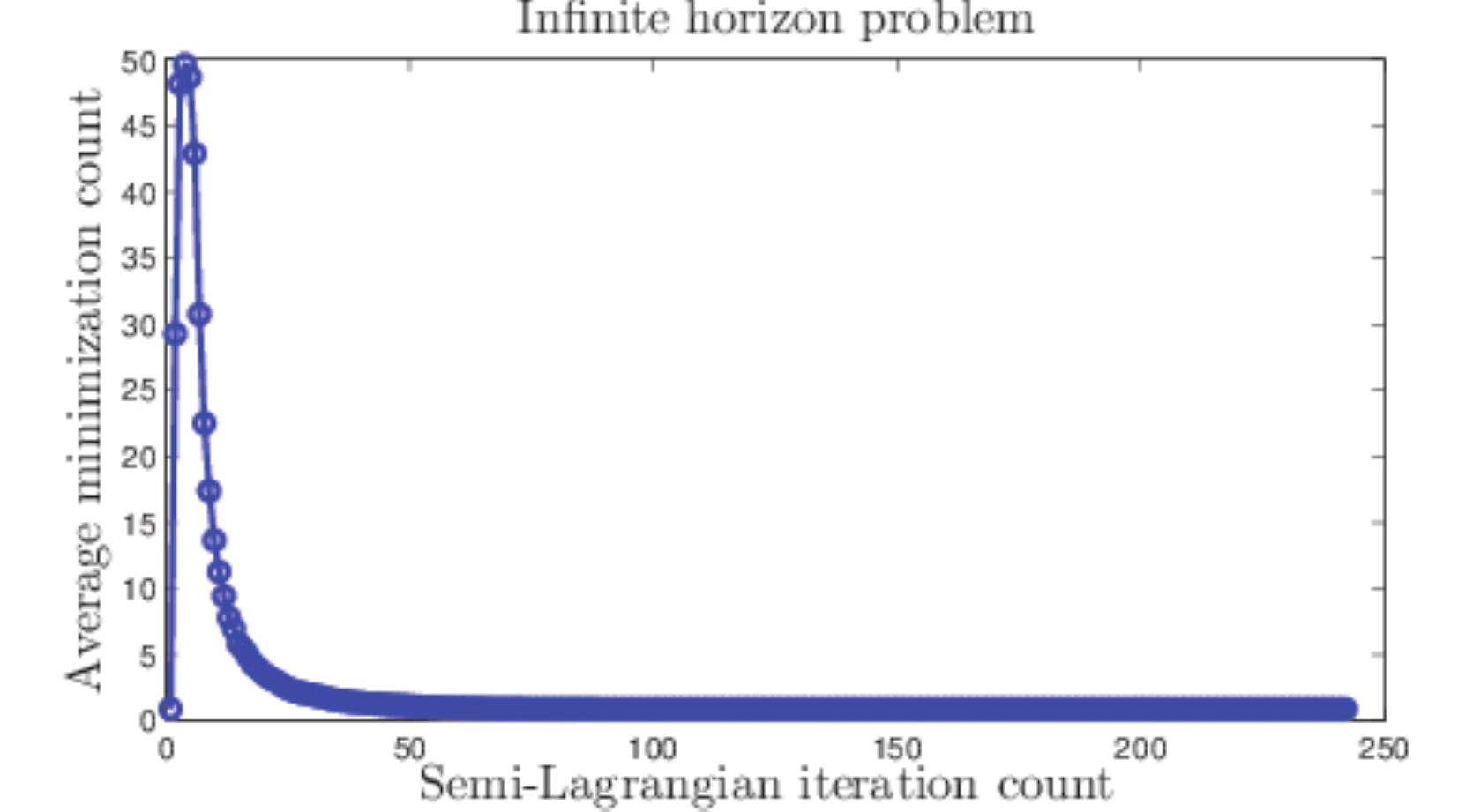,width=.45\linewidth,clip=}
  \caption{Subiteration count for 2D control problems with eikonal dynamics}\label{figeikcount}
\end{figure}

Having embedded the minimization routines within a semi-Lagrangian scheme, 
we show in Fig. \ref{figeikcount}  the evolution of the average iteration count per gridpoint (at
every fixed point iteration of the semi-Lagrangrian scheme), for both a minimum time and an infinite horizon optimal
control problem subject to eikonal dynamics. Note the
difference in the evolution of the subiteration count depending on the considered control
problem. In the minimum time problem nodes that have not received information propagating from the
optimality front are still minimized with irrelevant information until they are reached by the
optimality front, whereas in the infinite horizon case ``correct'' information is available for every node from the first iteration due to the presence of a running cost. In this latter case, the impact of the available information from the previous semi-Lagrangian iteration is similar to a warm start of the optimization routine.

\subsection{Infinite horizon problems with  $\ell_2$ running cost}
We present three different numerical tests for infinite horizon optimal control problem with cost functional and running cost given by
 \begin{eqnarray*}
	J(u,x)=\int_0^{\infty} l(x(s),u(s))e^{-\lambda s}ds\,,\quad\text{and }\quad
	l(u,x)=\frac12\|x\|_2^2+\frac{\gamma_2}{2}\|u\|_2^2\,,\quad\gamma_2>0,\lambda>0\,.
\end{eqnarray*}
 As common setting, the fixed point iteration is solved until
\[\|V^{n+1}-V^{n}\|\leq \frac15k^2\,,\quad n=1,2,3,\dots\,,\]
where $k$ stands for the space discretization parameter and the stopping tolerance for the inner optimization routine is set to $10^{-4}$. Further common
parameters are $\lambda=0.1$ and $\gamma_2=2$, and specific settings for every problem can be found
in Table \ref{part1}.
\begin{table}[htb]
\centering
  \begin{tabular}{cclr}
    \hline\\
    Test & $\Omega$ &\multicolumn{1}{c}{$U$}& \multicolumn{1}{c}{$h$}  \\
    \cmidrule(r){1-1} \cmidrule(lr){2-2} \cmidrule(lr){3-3} \cmidrule(lr){4-4}\\
    Test 1     &$[-1,1]^2$ & $\|u\|_2\leq 1$ & $\frac{\sqrt{2}}{4} k$ \\
    Test 2     &$[-1,1]^3$ & $\|u\|_2\leq 1$ & $\frac{1}{2}k$  \\
    Test 3     &$[-1,1]^3$ & $\|u\|_{\infty}\leq 0.3$ & $\frac{1}{5} k$  \\
    \hline 
  \end{tabular}
  \caption{Parameters for Tests 1-3.}\label{part1}
\end{table}

\subsubsection*{Test 1: 2D eikonal dynamics}
\begin{figure}[!ht]
  \centering
  \epsfig{file=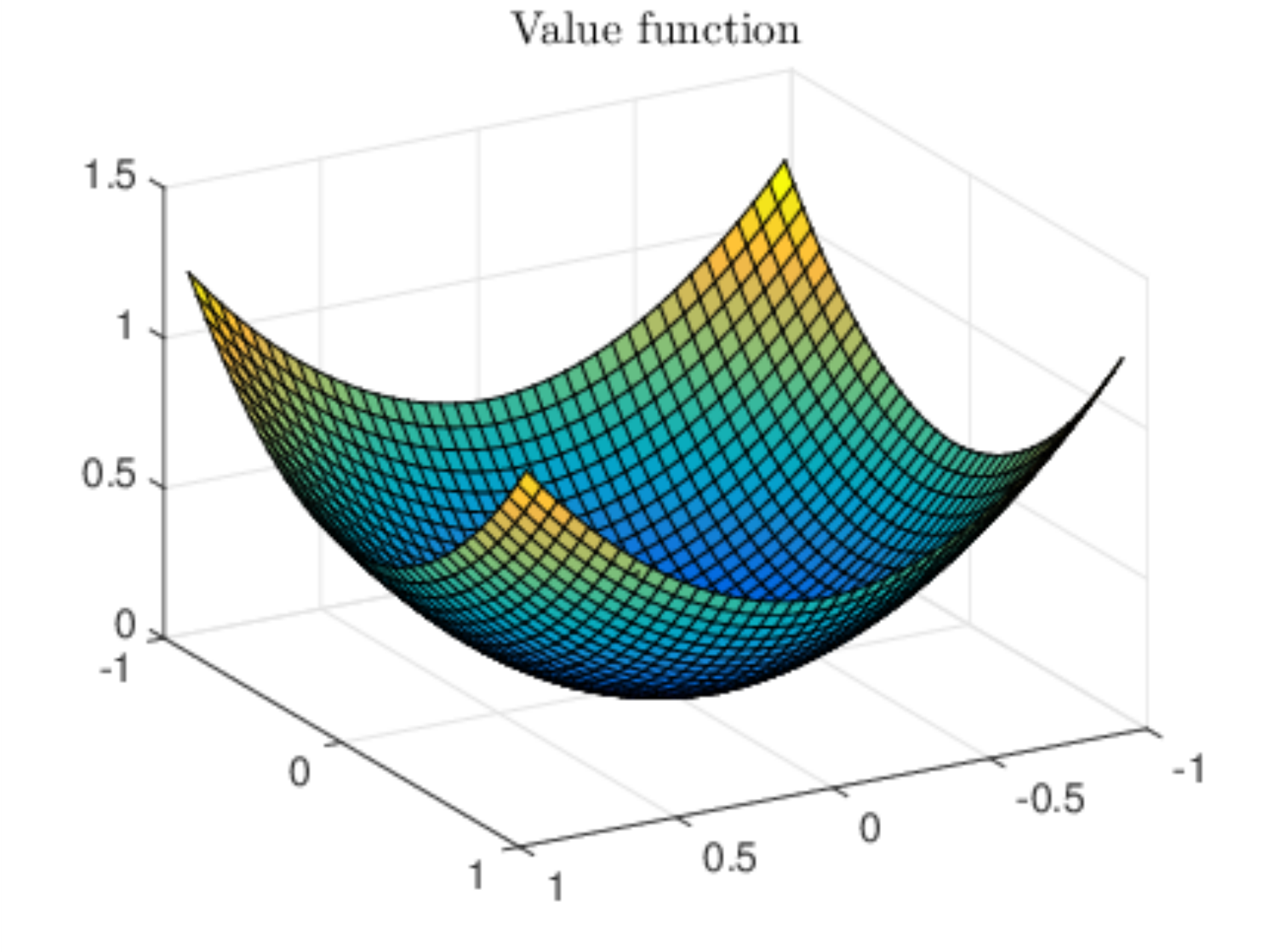,width=.45\linewidth,clip=}
    \epsfig{file=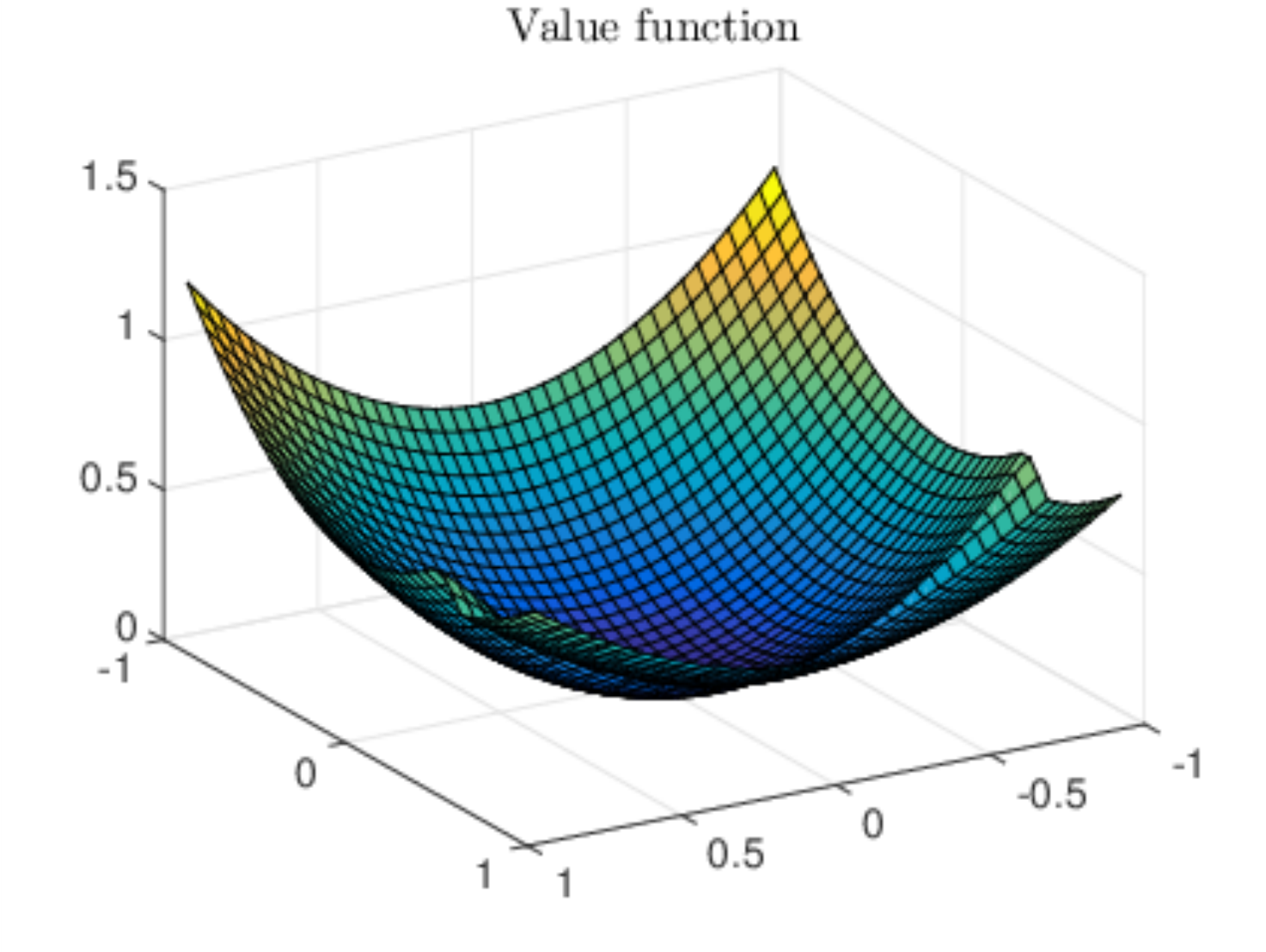,width=.45\linewidth,clip=}
     \epsfig{file=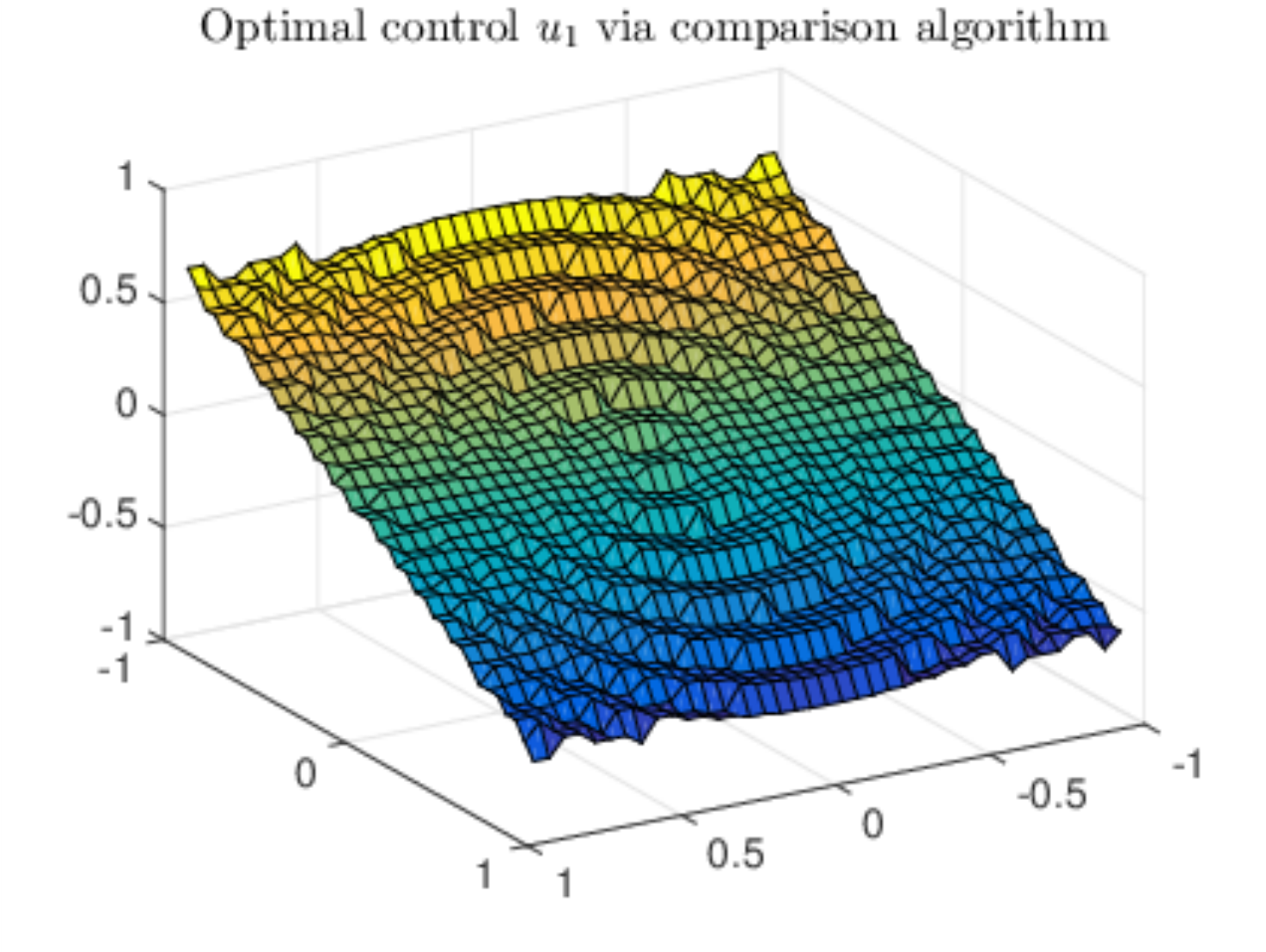,width=.45\linewidth,clip=}
    \epsfig{file=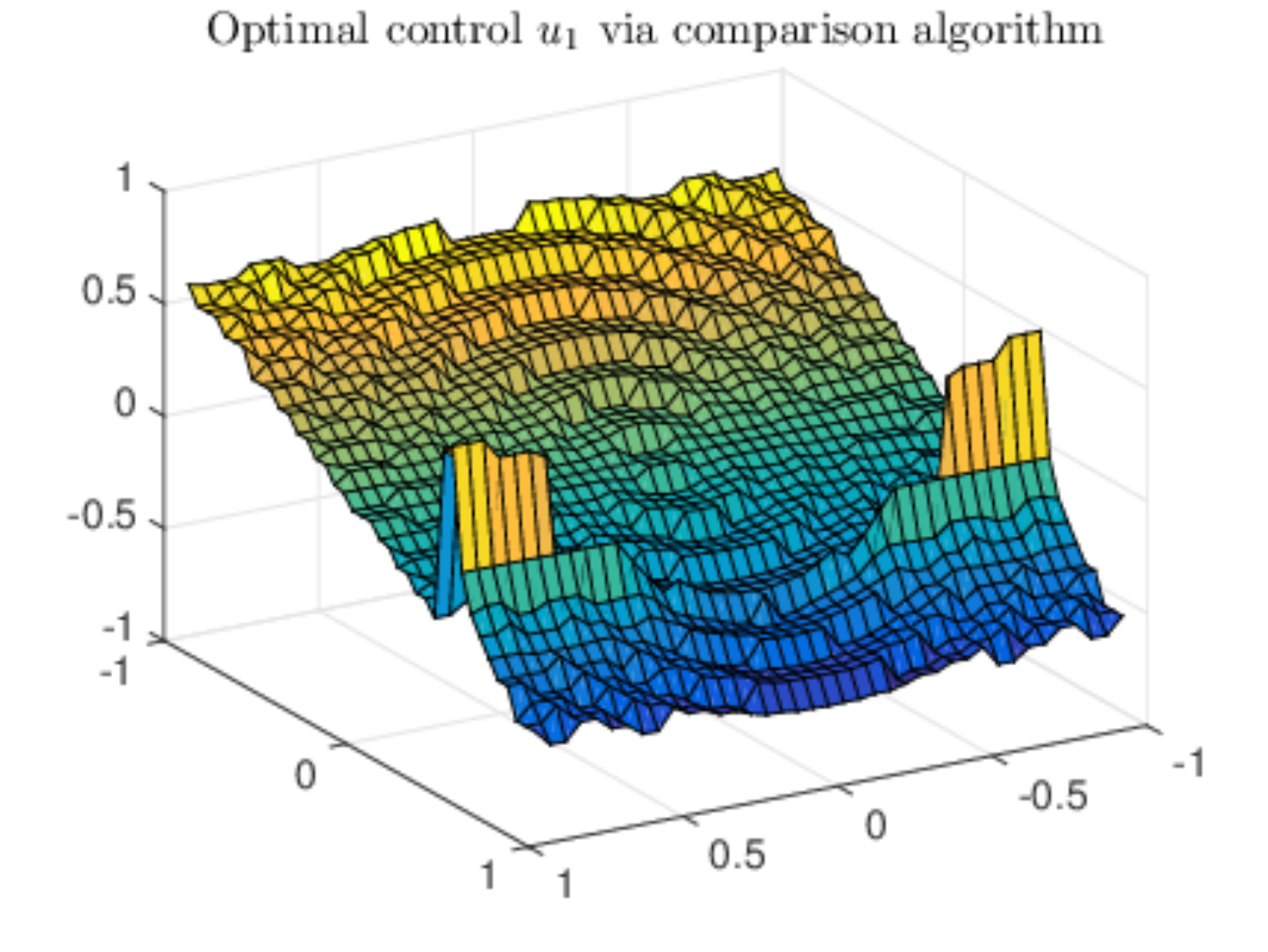,width=.45\linewidth,clip=}
  \epsfig{file=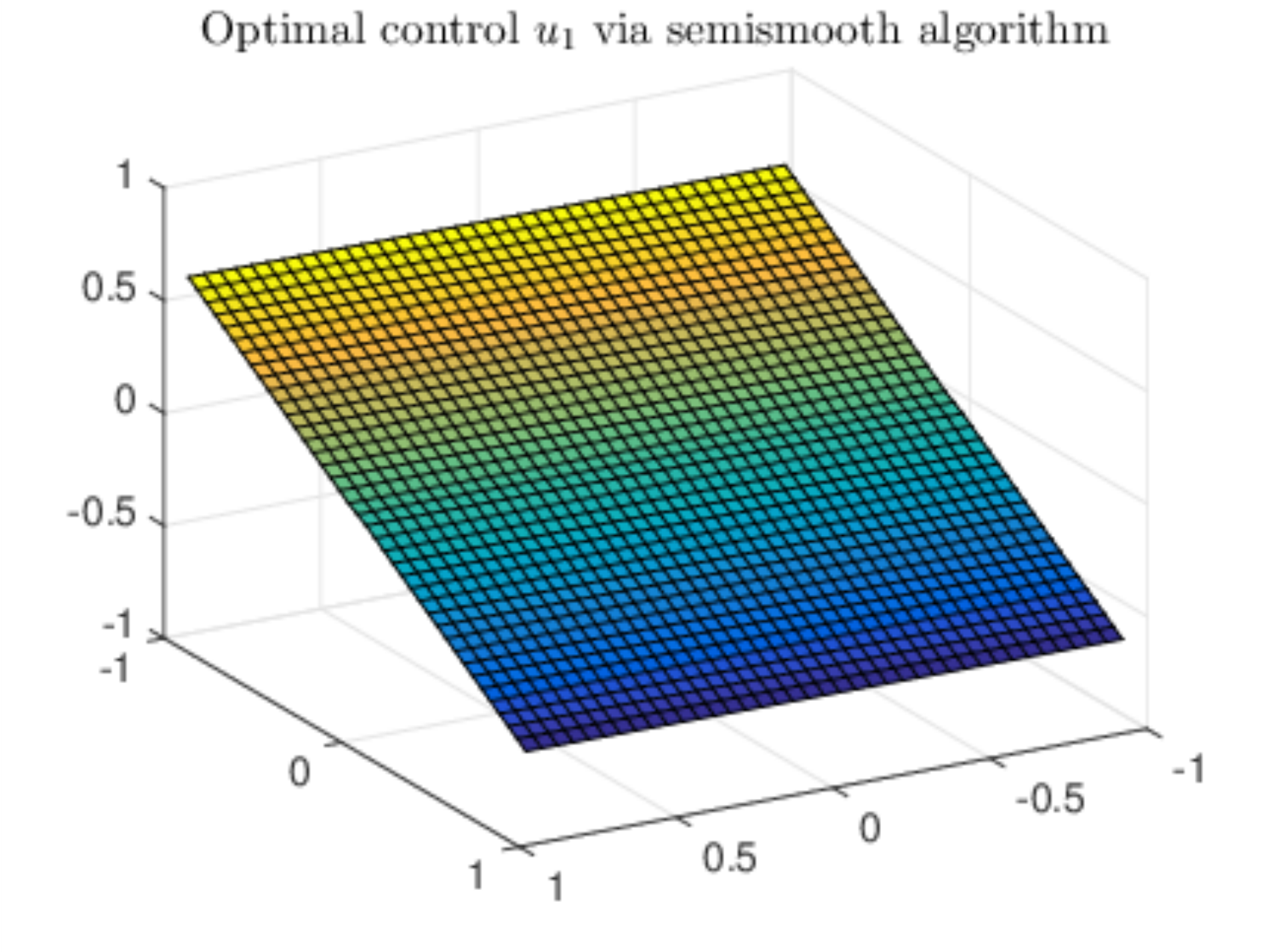,width=.45\linewidth,clip=}
  \epsfig{file=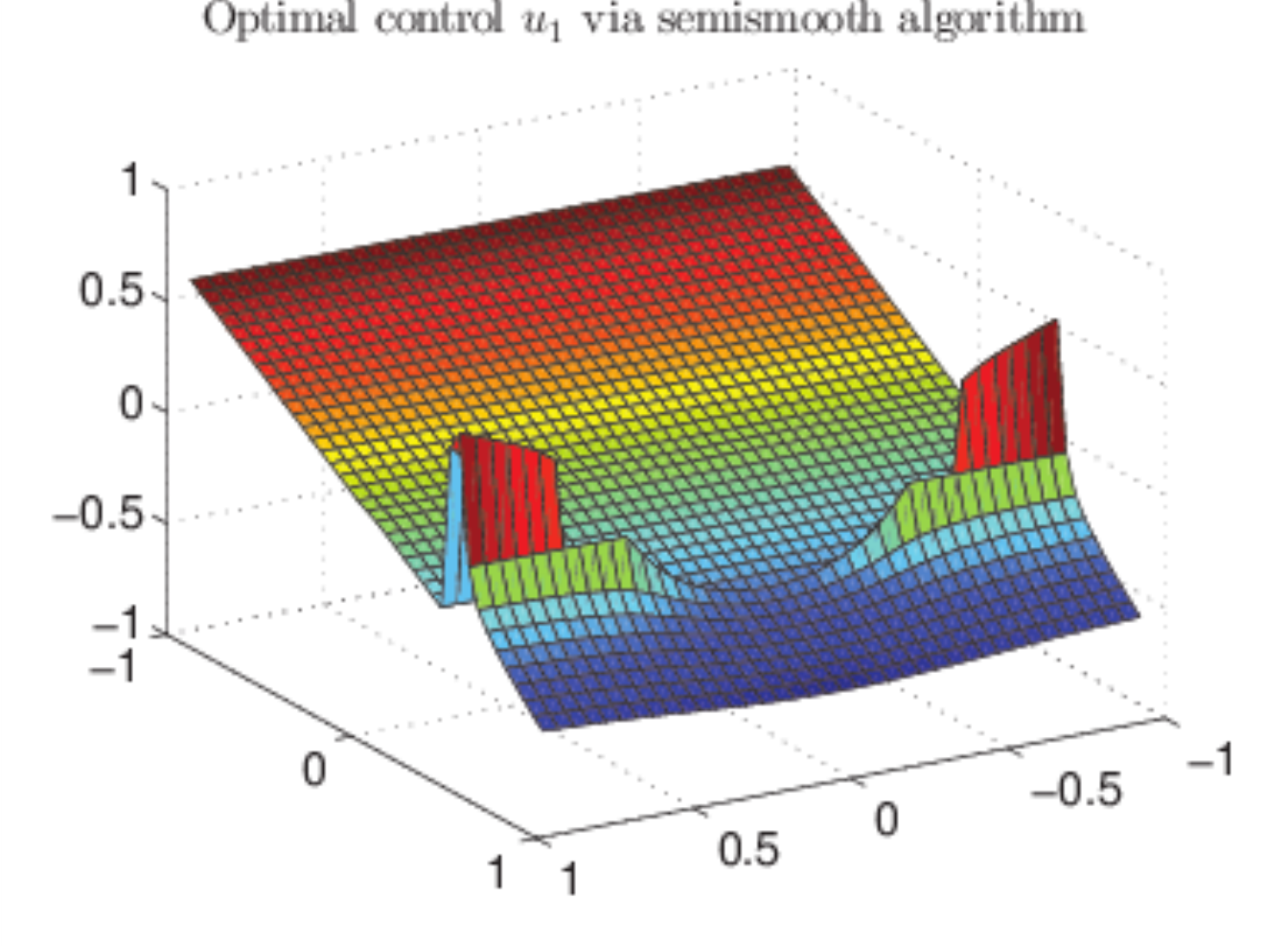,width=.45\linewidth,clip=}
  \caption{Infinite horizon control with 2D eikonal dynamics. Left: continuous dynamics. Right: discontinuous dynamics.}\label{test1}
\end{figure}

In this test we consider eikonal dynamics of the form
\begin{equation}
f(x,u)=(u_1,u_2)^T\,,\quad\|(u_1,u_2)\|_2\leq 1\,.
\end{equation}
We study the accuracy and performance of the semi-Lagrangian scheme with different minimization
routines: discretization of the control set and minimization by comparison, a semismooth Newton
method given in Algorithm \ref{ssn1}, and the approach by Chambolle and Pock given in Algorithm
\ref{cp2}. In order to make a fair study of the different routines, we choose the discrete set of
controls for the comparison algorithm such that the CPU time for the semismooth Newton method and
the Chambolle-Pock algorithm is almost the same. Table \ref{tabtest1} shows errors in both the value function and in the optimal control between the exact
solutions $v$ and $u$ and their
numerical approximations $V_h$ and $U_h$ for the different
schemes. Errors
are computed with respect to the exact solution of the Hamilton-Jacobi-Bellman equation, which is
not readily available in the literature. Since it is useful for numerical investigations, it is
provided in Appendix B.
The results show that we have similar CPU time of the approaches and independently of the meshsize,
the optimization-based schemes yield more accurate approximations of the value function and the
associated optimal control than the approach based on comparison. Further results are shown in Fig. \ref{test1}, where we also consider  nonhomogeneous eikonal dynamics 
\begin{equation}
f(x,u)=(1+\chi_{x_2>0.5}(x))(u_1,u_2)^T\,,\quad\|(u_1,u_2)\|_2\leq 1
\end{equation}
with $\chi_{x_2>0.5}(x)$ corresponds to the indicator function of the set $\{x=(x_1,x_2)\,|
x_2>0.5\}$. The figure shows that both approaches, the SL-scheme with a semismooth inner
optimization block and the one with a comparison-based routine, lead to very similar value functions. This is clearly not the case for the optimal control fields depicted in rows 2 and 3 of Fig. \ref{test1}. Even by a post-processing step it would be difficult to obtain the results in the third row from those in the second row.

\begin{center}
\begin{table}[htb]
\scalebox{0.84}{
\setlength{\tabcolsep}{1mm}
  \begin{tabular}{cp{1.7cm}p{1.9cm}p{1.7cm}p{1.7cm}p{1.9cm}p{1.7cm}}
  \hline\\
  & \multicolumn{3}{c}{$k=0.05$, \quad ($40^2$ DoF)}&  \multicolumn{3}{c}{$k=0.025$, \quad ($80^2$ DoF)} \\
\cmidrule(lr){2-4} \cmidrule(lr){5-7}\\
    Algorithm & CPU time&  $\|v-V_h\|_1$ & $\|u-U_h\|_1$ & CPU time&  $\|v-V_h\|_1$ & $\|u-U_h\|_1$\\
   \cmidrule(lr){1-1}  \cmidrule(lr){2-4} \cmidrule(lr){5-7}\\
    Comparison   &63.52 [s] &    3.12E-2 &   3.84E-2 &5.76E2 [s]&    1.96E-2 &   1.74E-2 \\
    Semismooth Newton &77.25 [s]&    2.62E-2 &  1.61E-2   &7.27E2 [s] &    1.36E-2 &  7.21E-3\\
    Chambolle-Pock&63.05 [s] & 2.60E-2 &   1.42E-2 &5.77E2 [s] &   1.36E-2 &  6.83E-3\\
    \hline 
  \end{tabular}}
  \caption{Infinite horizon control of 2D eikonal dynamics. CPU time and errors for a semi-Lagrangian
  scheme with different minimization routines. The comparison algorithm was run with a discrete set
  of 1280 control points in every node.}\label{tabtest1}
\end{table}
\end{center}

\subsubsection*{Test 2: 3D eikonal dynamics}
In order to illustrate that our approach can be extended to higher dimensions, we consider the infinite horizon optimal control of three dimensional eikonal dynamics
$$f(x,u)=(u_1,u_2,u_3)^T\,,\quad\|(u_1,u_2,u_3)\|_2\leq 1\,.$$
Errors and CPU time are shown in Table \ref{tab20}. Since the values for the semismooth Newton and the Chambolle-Pock algorithms are similar, we only include the data for this latter one. With the Chambolle-Pock algorithm we obtain more accurate solutions of the control field for a similar amount of CPU time as the comparison approach.

\begin{center}
\begin{table}[htb]
\scalebox{.88}{
\setlength{\tabcolsep}{1mm}
  \begin{tabular}{cp{1.7cm}p{1.9cm}p{1.7cm}p{1.7cm}p{1.9cm}p{1.7cm}}
  \hline\\
     & \multicolumn{3}{c}{$k=0.1$, \quad ($20^3$ DoF)} & \multicolumn{3}{c}{$k=0.05$, \quad ($40^3$ DoF)} \\
\cmidrule(lr){2-4} \cmidrule(lr){5-7}\\
    Algorithm & CPU time&  $\|v-V_h\|_1$ & $\|u-U_h\|_1$ & CPU time&  $\|v-V_h\|_1$ & $\|u-U_h\|_1$\\
   \cmidrule(lr){1-1}  \cmidrule(lr){2-4} \cmidrule(lr){5-7}\\
    Comparison   &93.76 [s] &    2.49E-2 &   2.46E-2 &2.89E2 [s]&    1.02E-2 &   1.89E-2 \\
    Chambolle-Pock&97.25 [s] & 9.92E-3 &   2.07E-2 &2.01E2 [s] &   5.66E-3 &  1.22E-2\\
    \hline 
  \end{tabular}}
  \caption{Infinite horizon control of 3D eikonal dynamics. CPU time and errors for a semi-Lagrangian scheme with different minimization routines. The comparison algorithm was run with a discrete set of 5120 control points.}\label{tab20}
\end{table}
\end{center}


\subsubsection*{Test 3: Triple integrator with two control variables}
\begin{figure}[!ht]
  \centering
  \epsfig{file=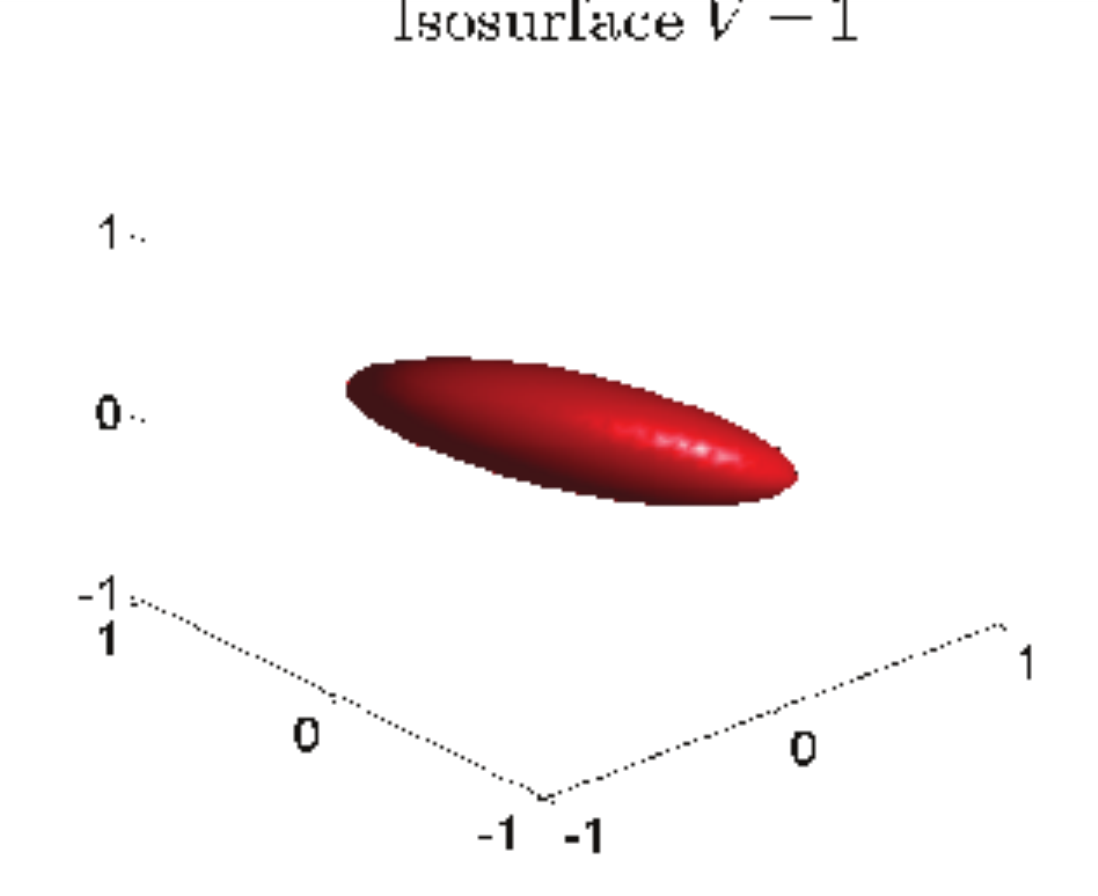,width=.3\linewidth,clip=}
  \epsfig{file=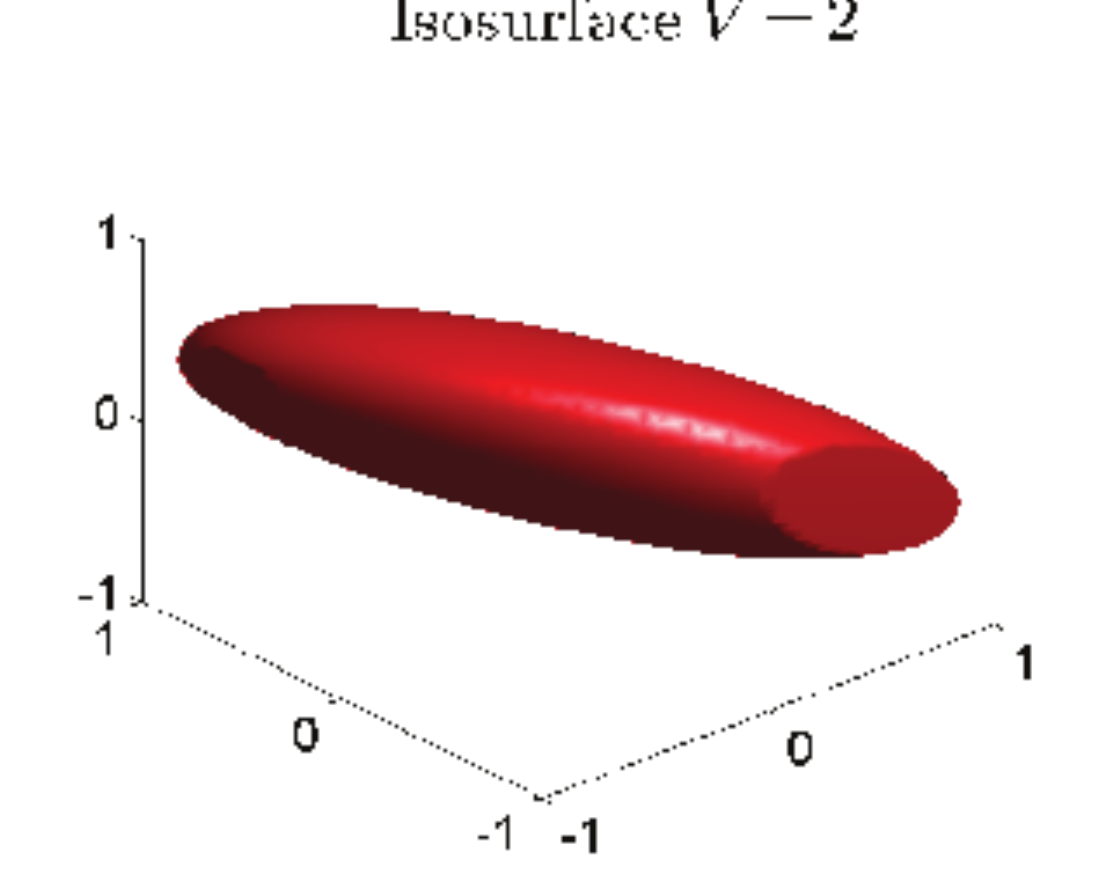,width=.3\linewidth,clip=}
  \epsfig{file=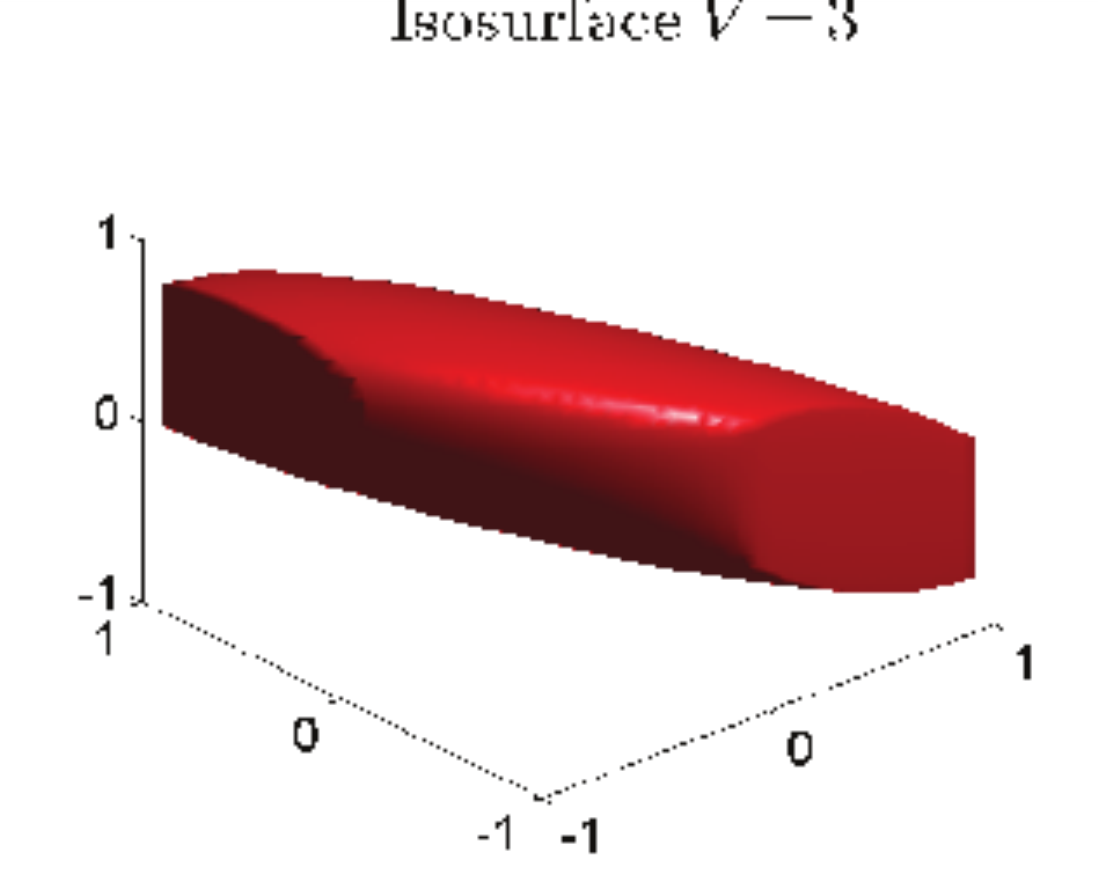,width=.3\linewidth,clip=}

  \epsfig{file=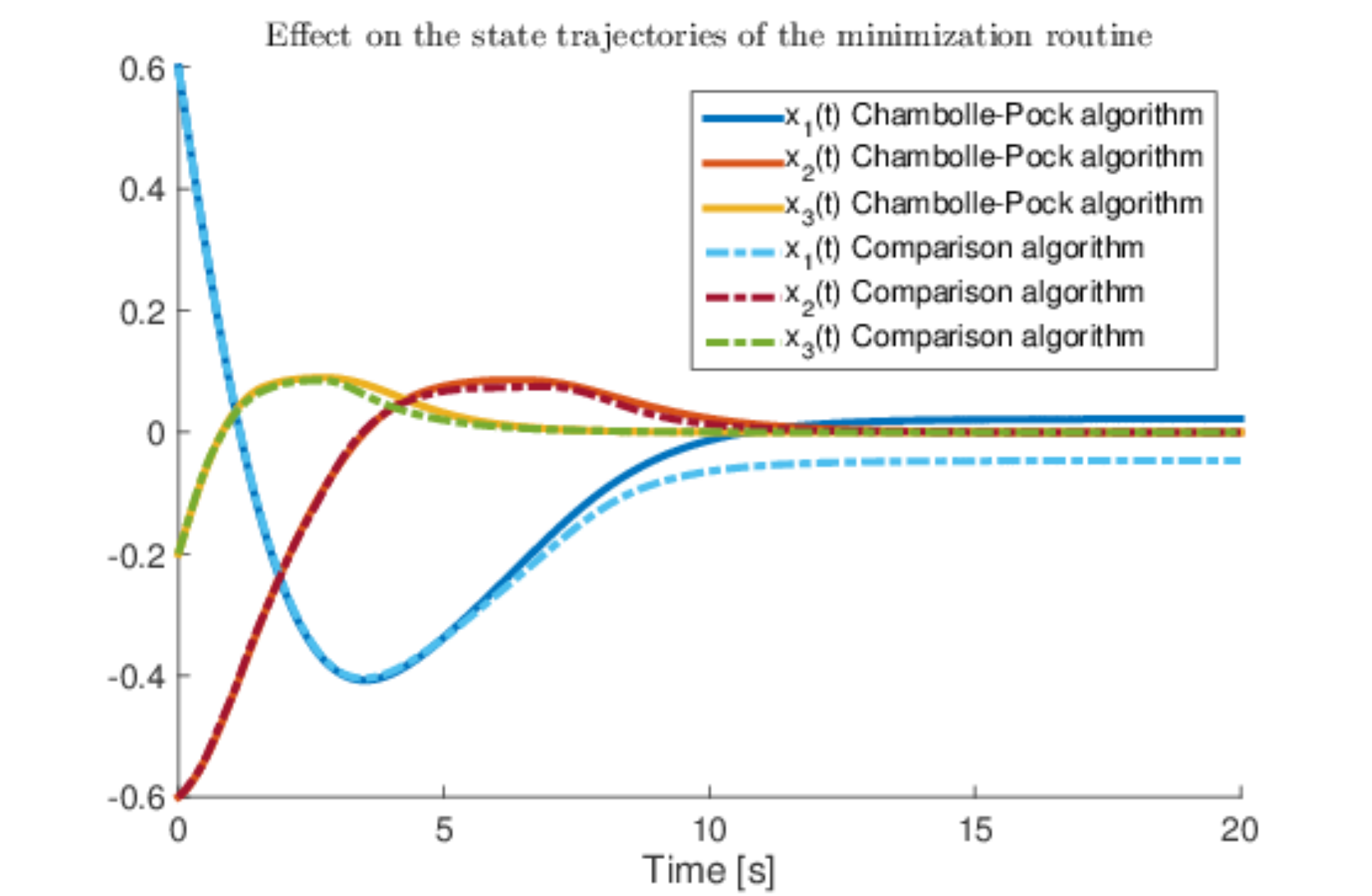,width=.45\linewidth,clip=}
  \epsfig{file=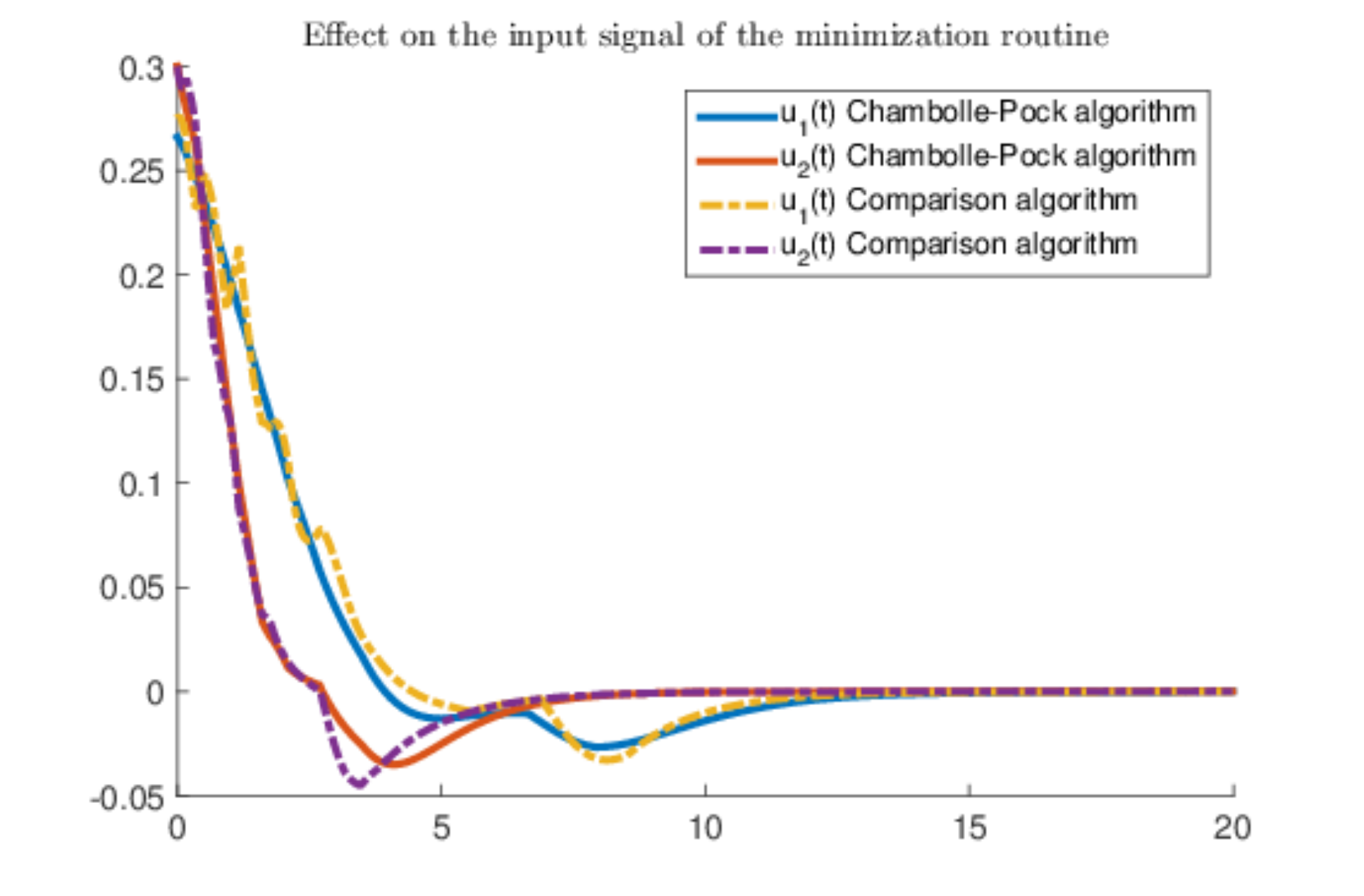,width=.45\linewidth,clip=}
  \epsfig{file=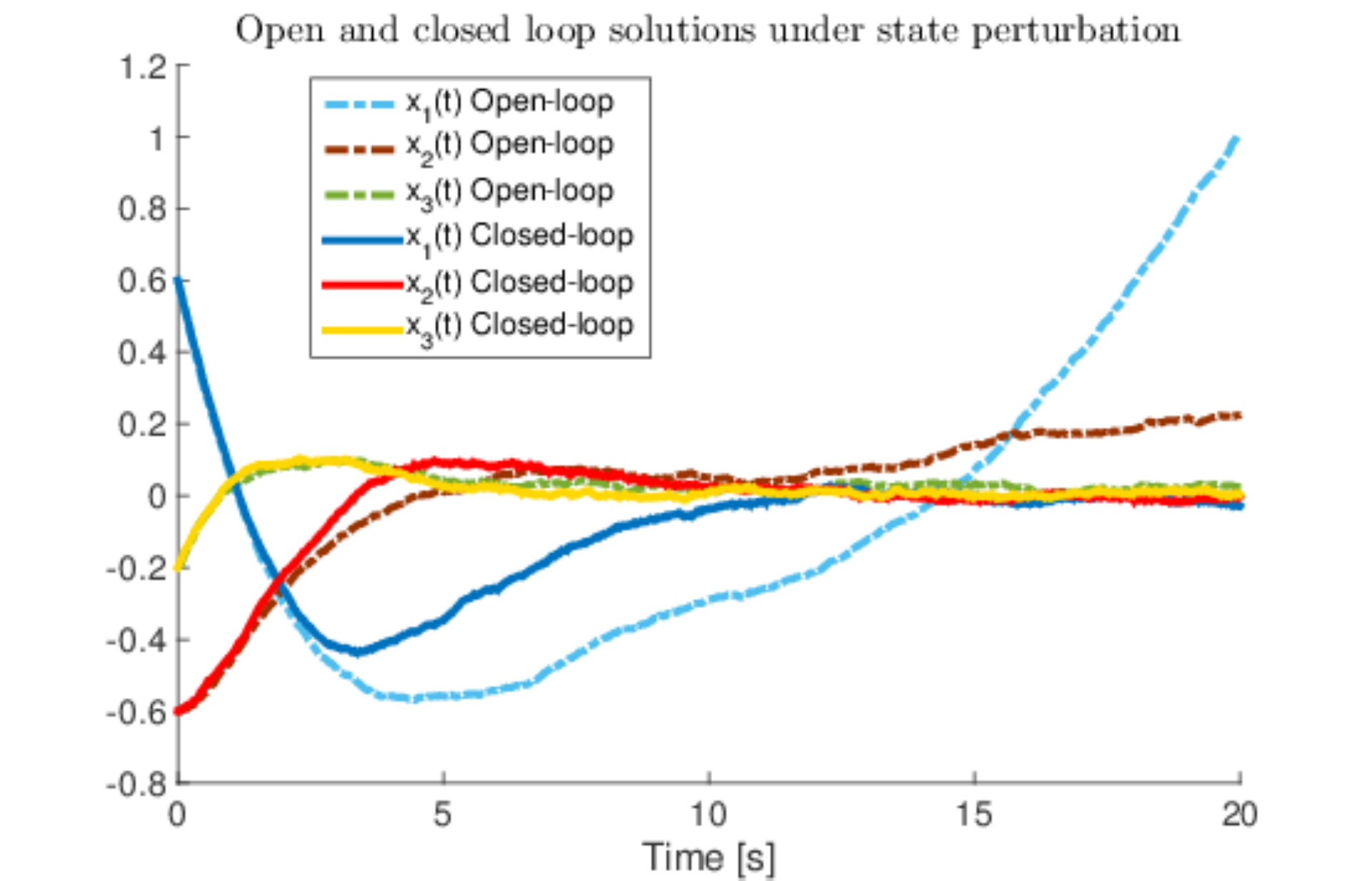,width=.45\linewidth,clip=}
  \epsfig{file=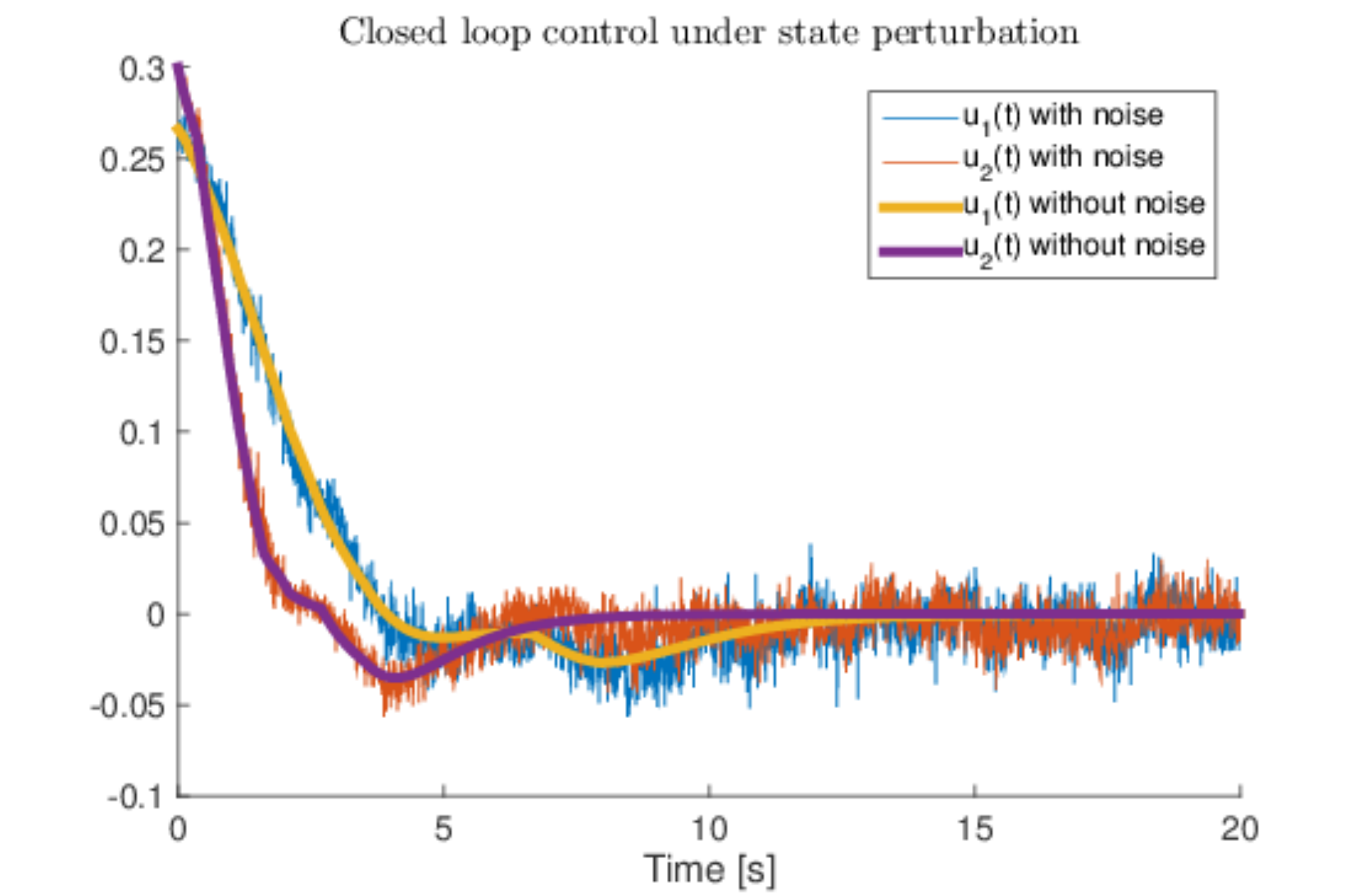,width=.45\linewidth,clip=}
  \caption{Triple integrator with two controls, numerical results with $k=0.05$. Top: different isosurfaces. Middle: differences in the optimal control lead to different trajectories. Bottom: the feedback approach leads to robustness with respect to noise in the dynamics.}\label{test3}
\end{figure}

In this test the dynamics are given by a triple integrator with two control variables subject to box constraints
$$f(x,u)=(x_2,x_3+u_1,u_2)^T\,,\quad|u_1|\leq a,\;|u_2|\leq b,\,\quad a,b\in \mathds{R}.$$
The purpose of this example is to stress that the minimization strategy that we have introduced can
be also applied to non-eikonal dynamics, where the correspondence between octants in the state space
and the control field is not trivial. For the sake of completeness the control space decomposition
procedure for this example can be found in Appendix A. In this particular case, every octant will
have associated a different rectangular sector in the control space for its arrival points. Results
for the value function, optimal controls and trajectories are shown in Fig. \ref{test3}. In the
second row of this figure, we observe distinct differences between the optimal controls obtained
from the Chambolle-Pock and comparison-based algorithms. These differences in the control lead to
different approximated steady states as it can be seen in the left of the second row of Fig. \ref{test3}. 
To highlight also the effect of closed-loop control we carried out an experiment with additive
structural and output noise. In the third row of Fig.~\ref{test3}, the resulting states from open and closed-loop control can be compared.

\subsection{Infinite horizon problems with  combined $\ell_2$ and $\ell_1$-costs}
We present two numerical tests for infinite horizon optimal control problems with cost function and running cost given by
 \begin{align*}
   &J(u,x)=\int_0^{\infty} l(x(s),u(s))e^{-\lambda s}ds\;,\\
	&l(u,x)=\frac12\|x\|_2^2+\frac{\gamma_2}{2}\|u\|_2^2+\gamma_1\|u\|_1\,,\quad\gamma_2>0,\;
        \gamma_1>0,\;\lambda>0\,.
\end{align*}
 The stopping rule for the fixed point iteration is defined in the same way as in the previous set
 of examples. All the optimization-based routines have been solved with a semismooth Newton method
 as the inner block, with a regularization parameter $\varepsilon=1e-3$. Further common parameters
 are $\lambda=0.1$ and $\gamma_2=2$, and specific settings for every problem can be found in Table \ref{part2}.
\begin{table}[htb]
\centering
  \begin{tabular}{lllcc}
    \hline\\
    Test & $\Omega$ &\multicolumn{1}{c}{$U$} & h & k  \\
    \cmidrule(r){1-1} \cmidrule(lr){2-2} \cmidrule(lr){3-3} \cmidrule(lr){4-4}  \cmidrule(lr){5-5}\\
    Test 4     &$[-1,1]^2$ & $\|u\|_2\leq 1$ & $\frac{\sqrt{2}}{4} k$ & 0.025\\
    Test 5     &$[0,2\pi]^3$ & $\|u\|_{\infty}\leq 0.3$ & $0.2 k$ & 0.2 \\
    \hline 
  \end{tabular}
  \caption{Parameters for Tests 4-5.}\label{part2}
\end{table}

\subsubsection*{Test 4: 2D eikonal dynamics}
\begin{figure}[!ht]
  \centering
  \epsfig{file=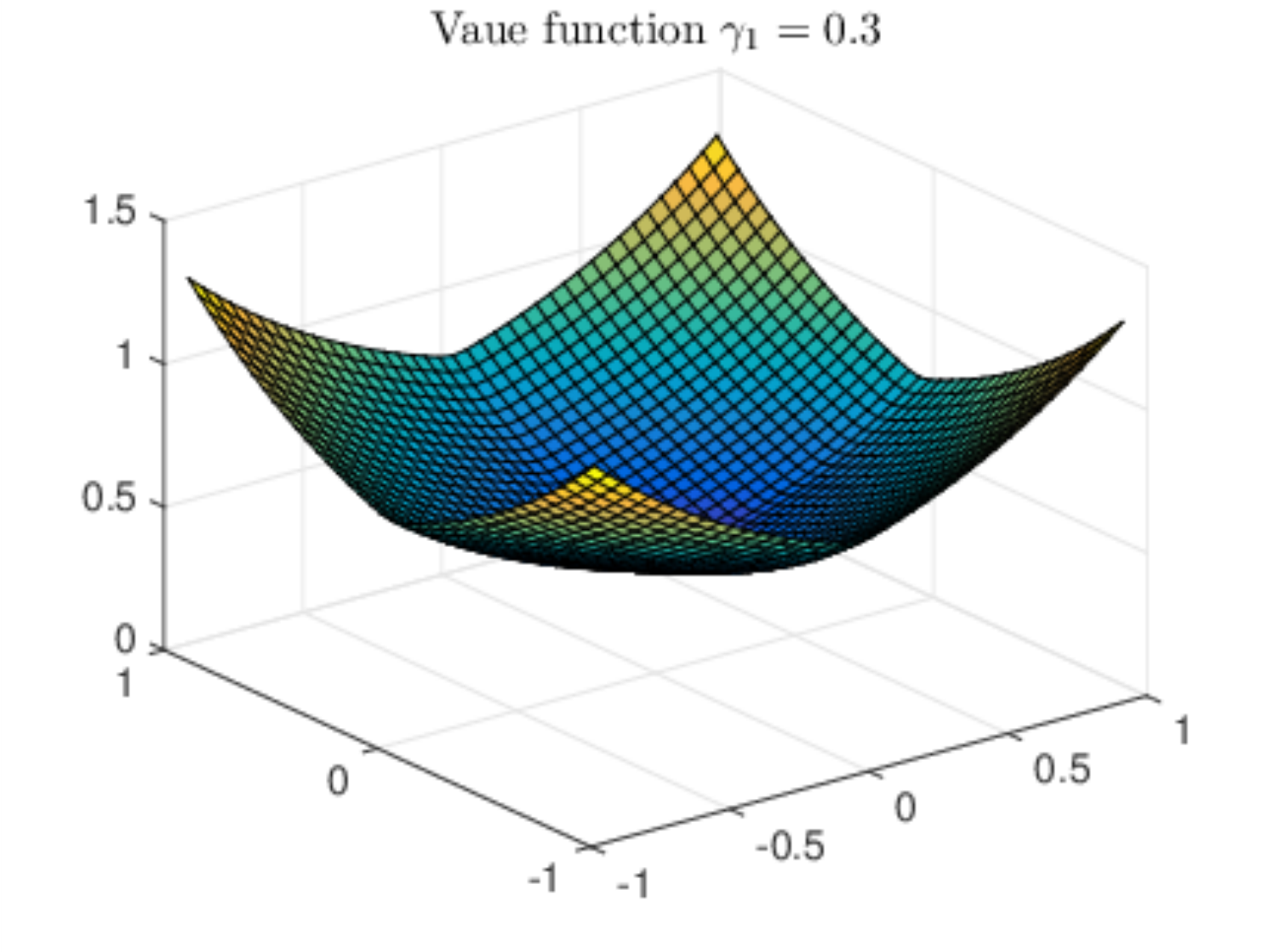,width=.48\linewidth,clip=}
  \epsfig{file=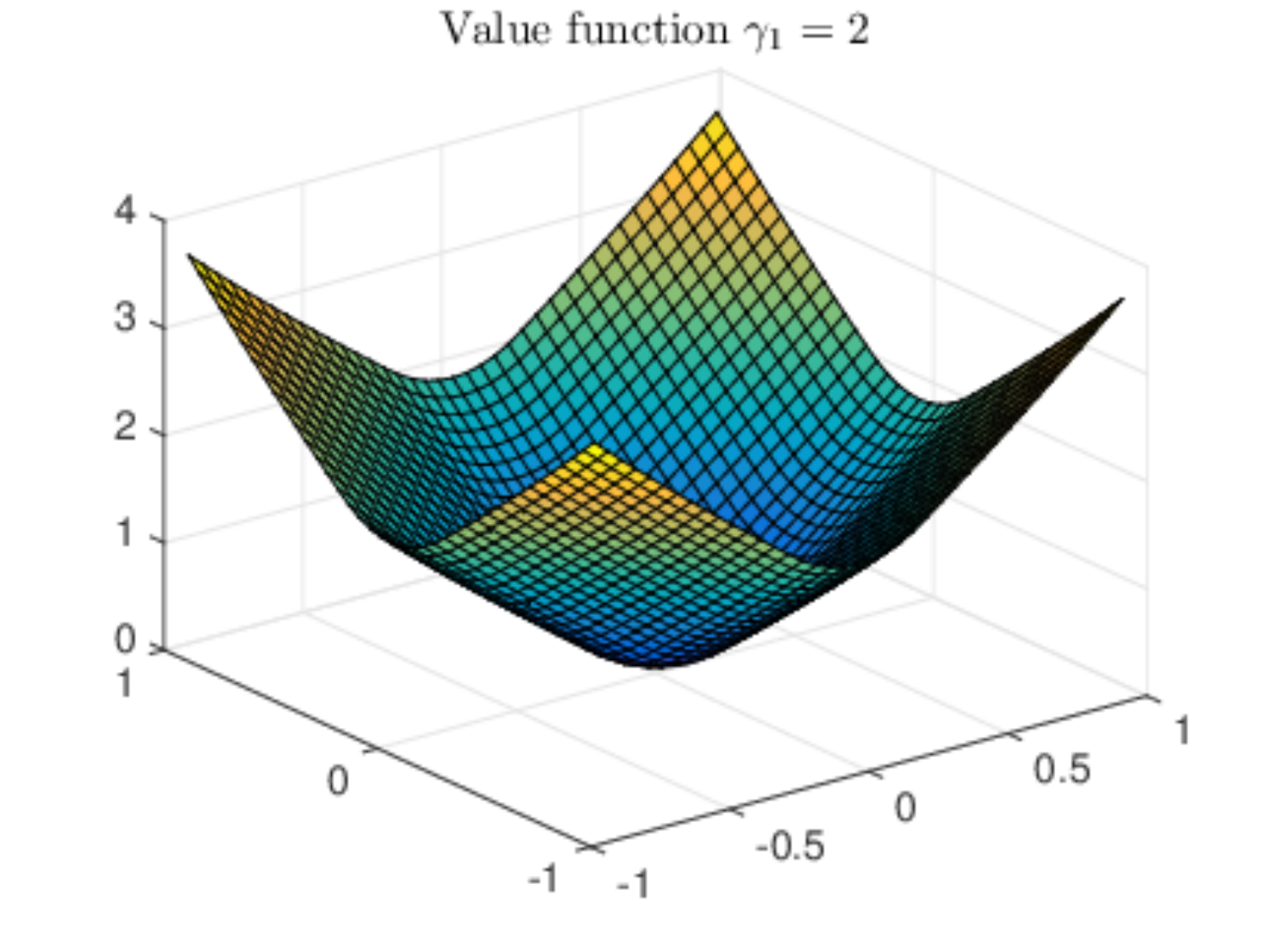,width=.48\linewidth,clip=}
  \epsfig{file=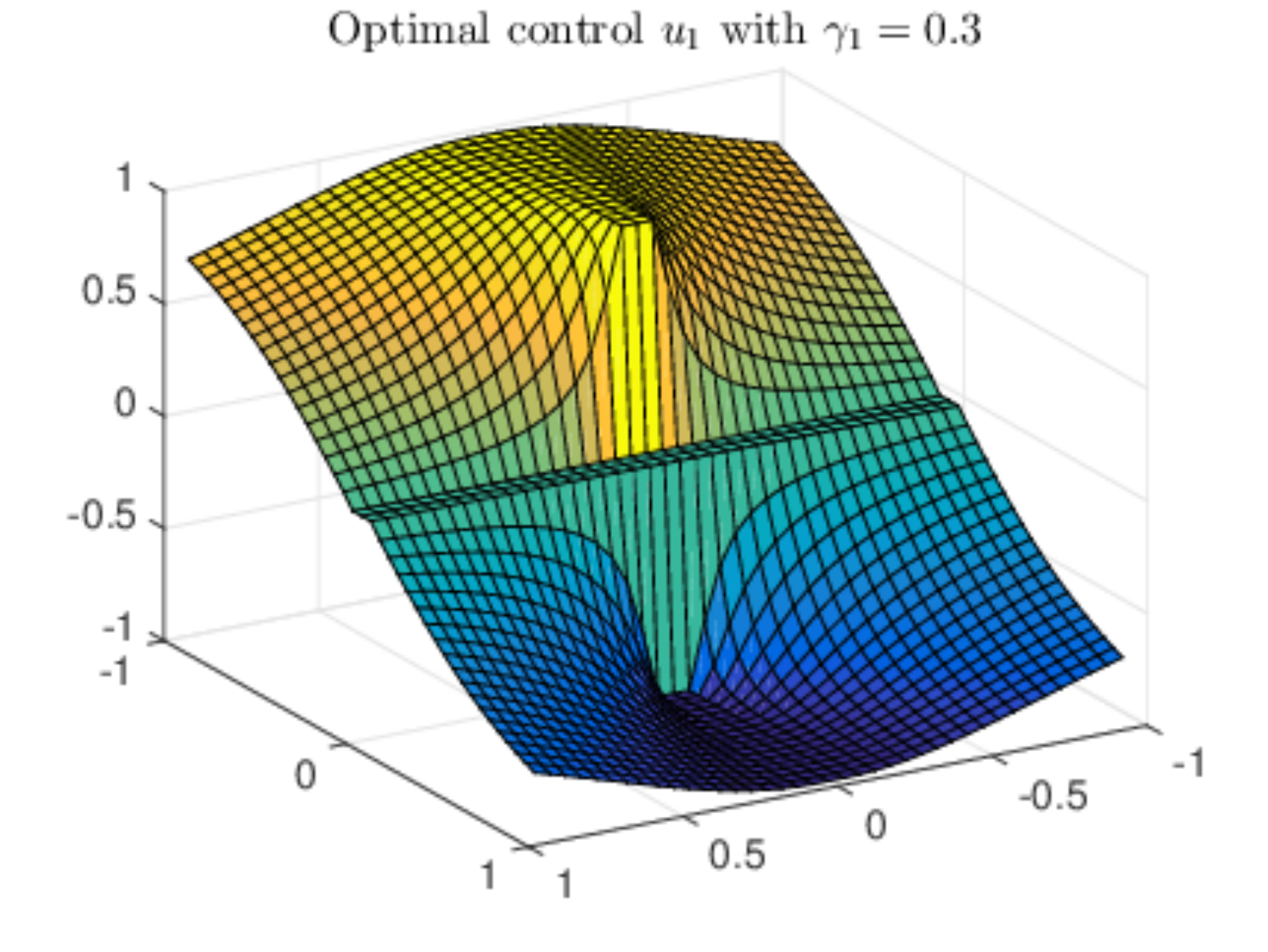,width=.48\linewidth,clip=}
  \epsfig{file=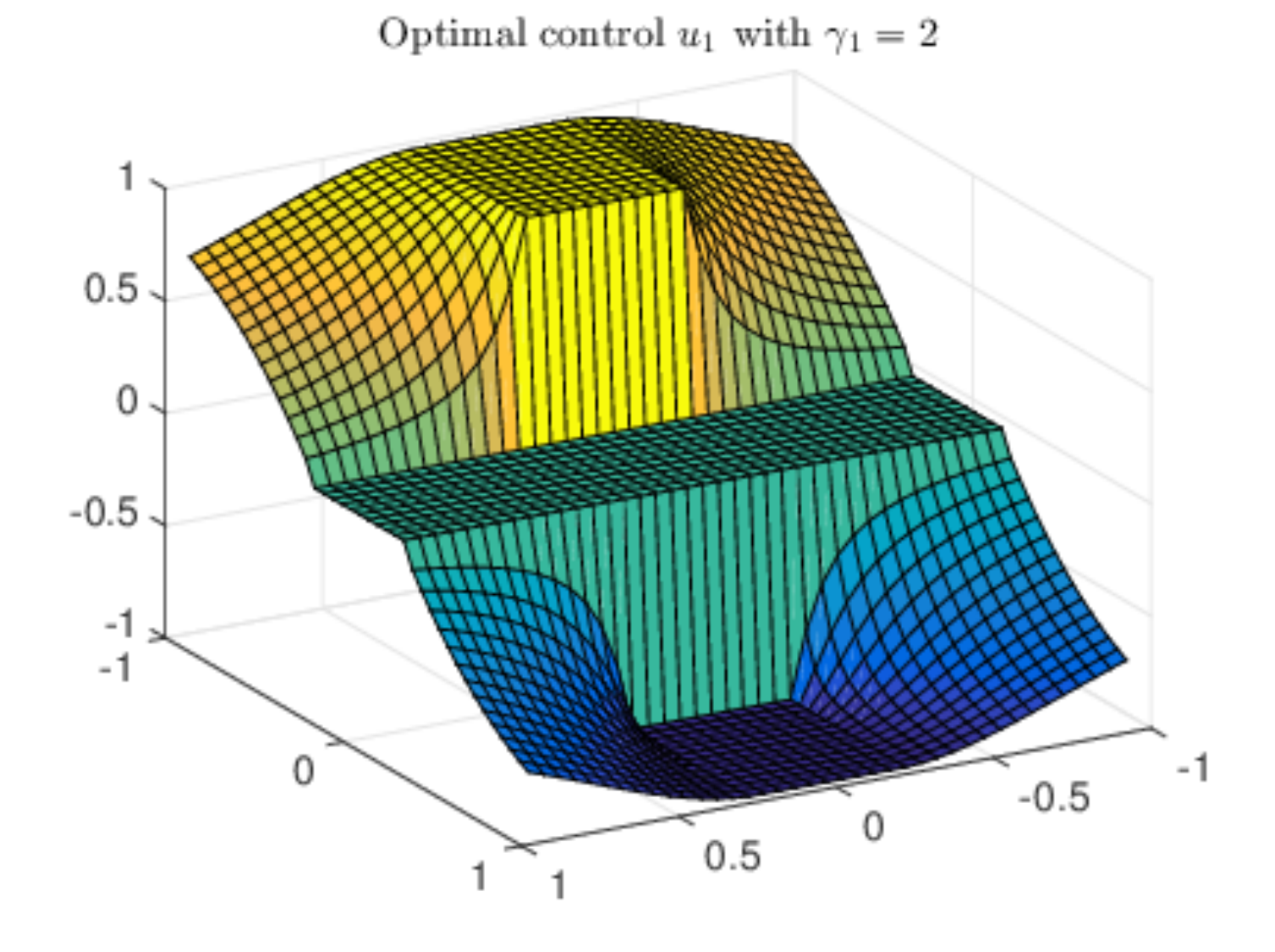,width=.48\linewidth,clip=}
  \caption{Sparse control of eikonal dynamics. Top: inclusion of an $\ell_1$-cost breaks the radially symmetric structure of the solution. Bottom: the sparsity of the optimal control translates into a zero band around of the origin, depending on the weight $\gamma_1$.}\label{test4}
\end{figure}

This test considers the same two dimensional dynamics as in Test 1, the only difference being the
inclusion of an $\ell_1$-term in the cost functional. For the inner optimization block we
apply the semismooth Newton method presented in Algorithm \ref{ssn2}. Results are shown in Fig.~\ref{test4}, where differences in the shape of the value
function can be observed as the $\gamma_1$ weight increases. In the second row of Fig. \ref{test4},
the effect of sparsity on the first control component can be seen from the fact that it is
identically zero in a band around the origin. Moreover, the width of the band increases with
$\gamma_1$. As for the value function, introducing the $\ell_1$-term breaks its radially symmetric
structure, as it is shown in the first row of Fig.~\ref{test4}.

\subsubsection*{Test 5: 3D car model}
\begin{figure}[!ht]
  \centering
   \epsfig{file=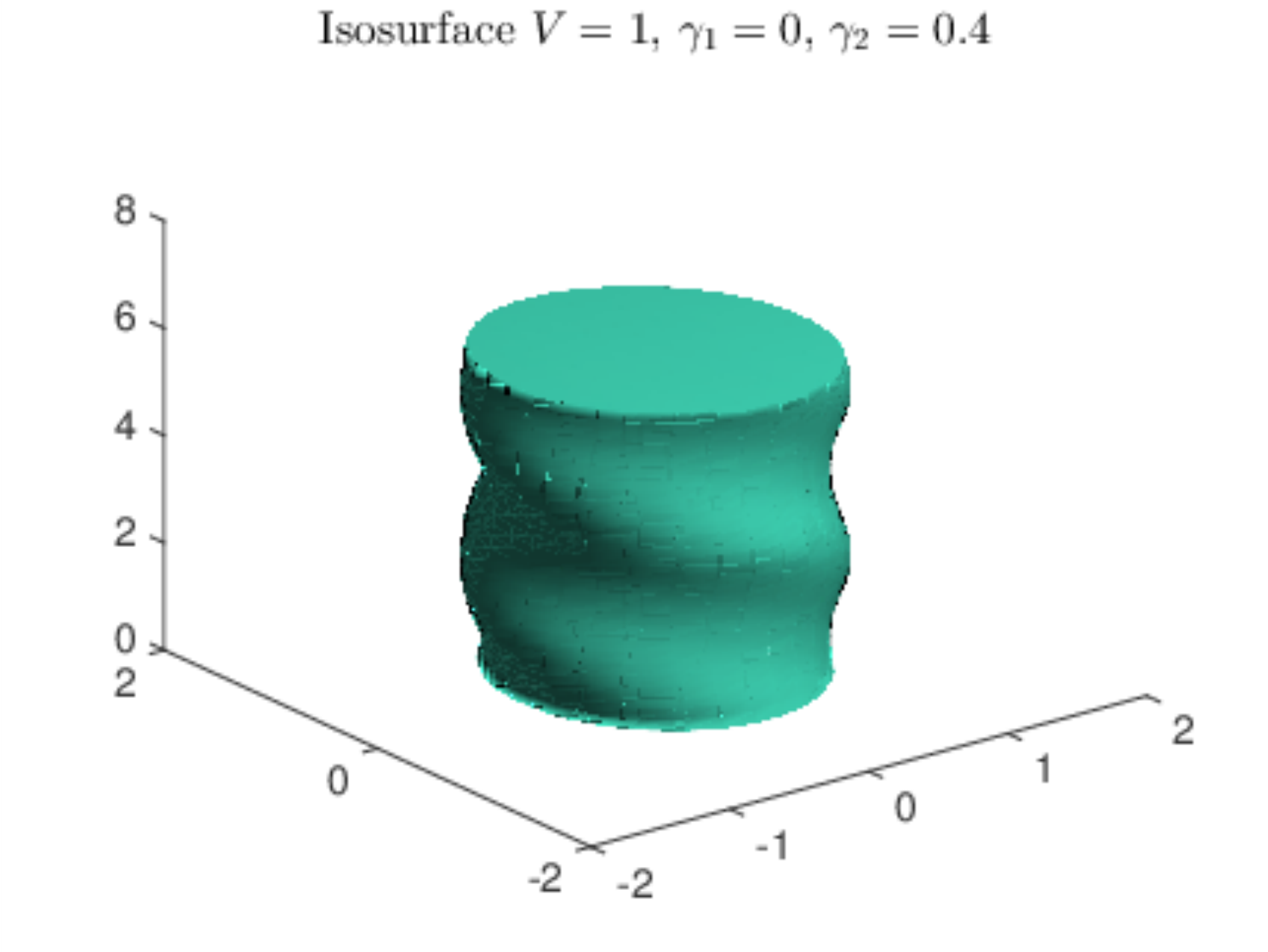,width=.48\linewidth,clip=}
  \epsfig{file=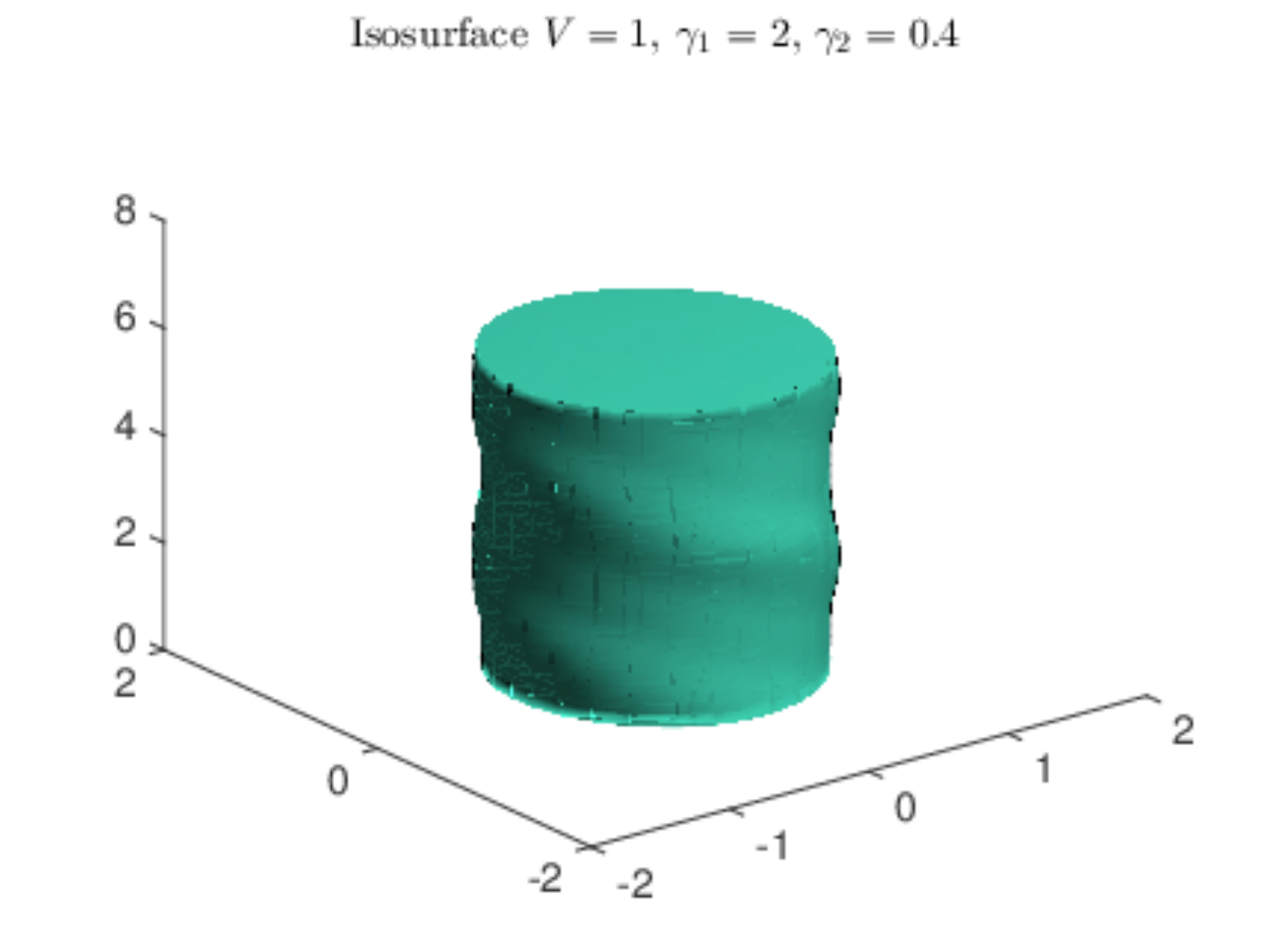,width=.48\linewidth,clip=}\\
  \epsfig{file=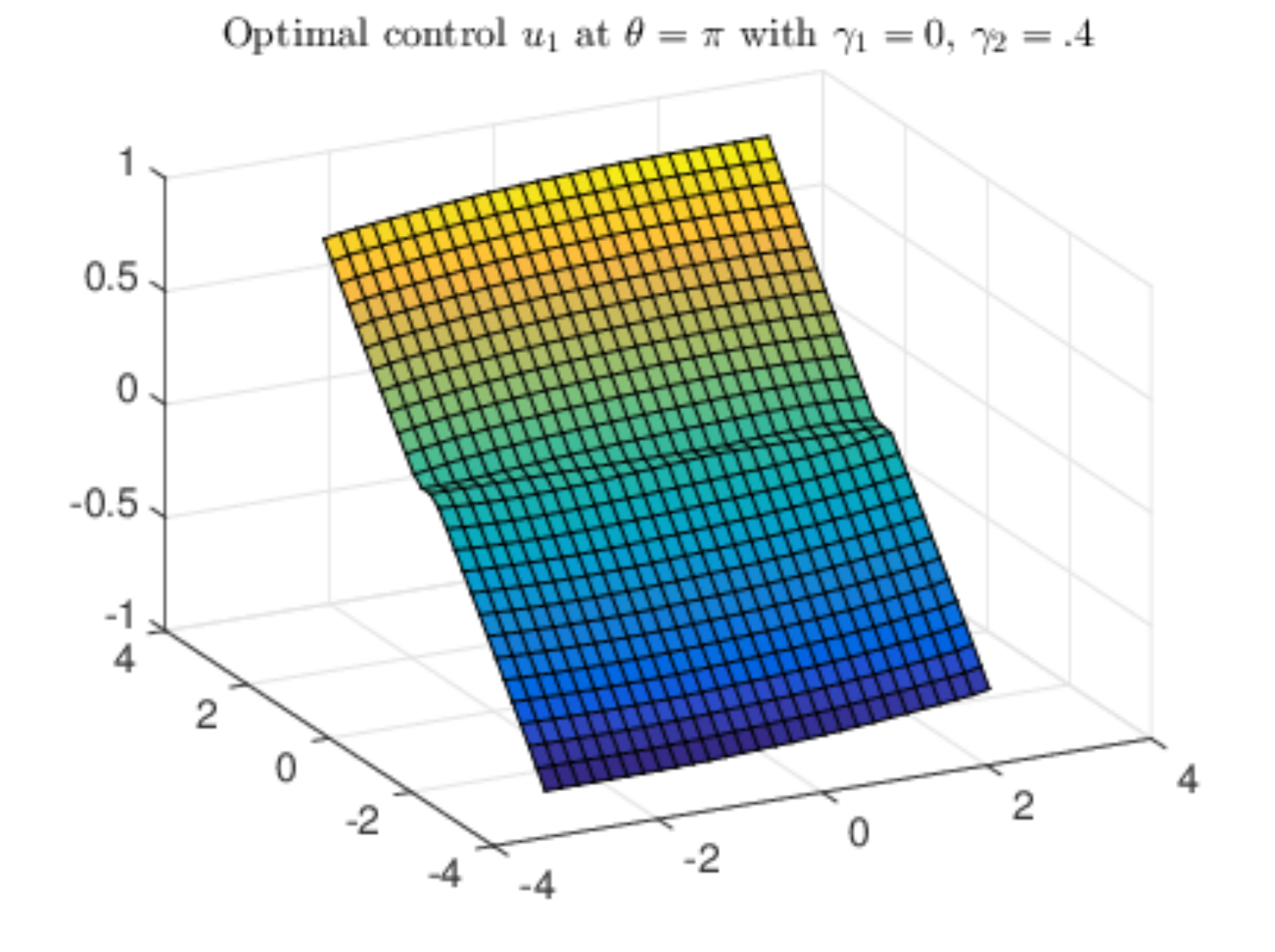,width=.48\linewidth,clip=}
  \epsfig{file=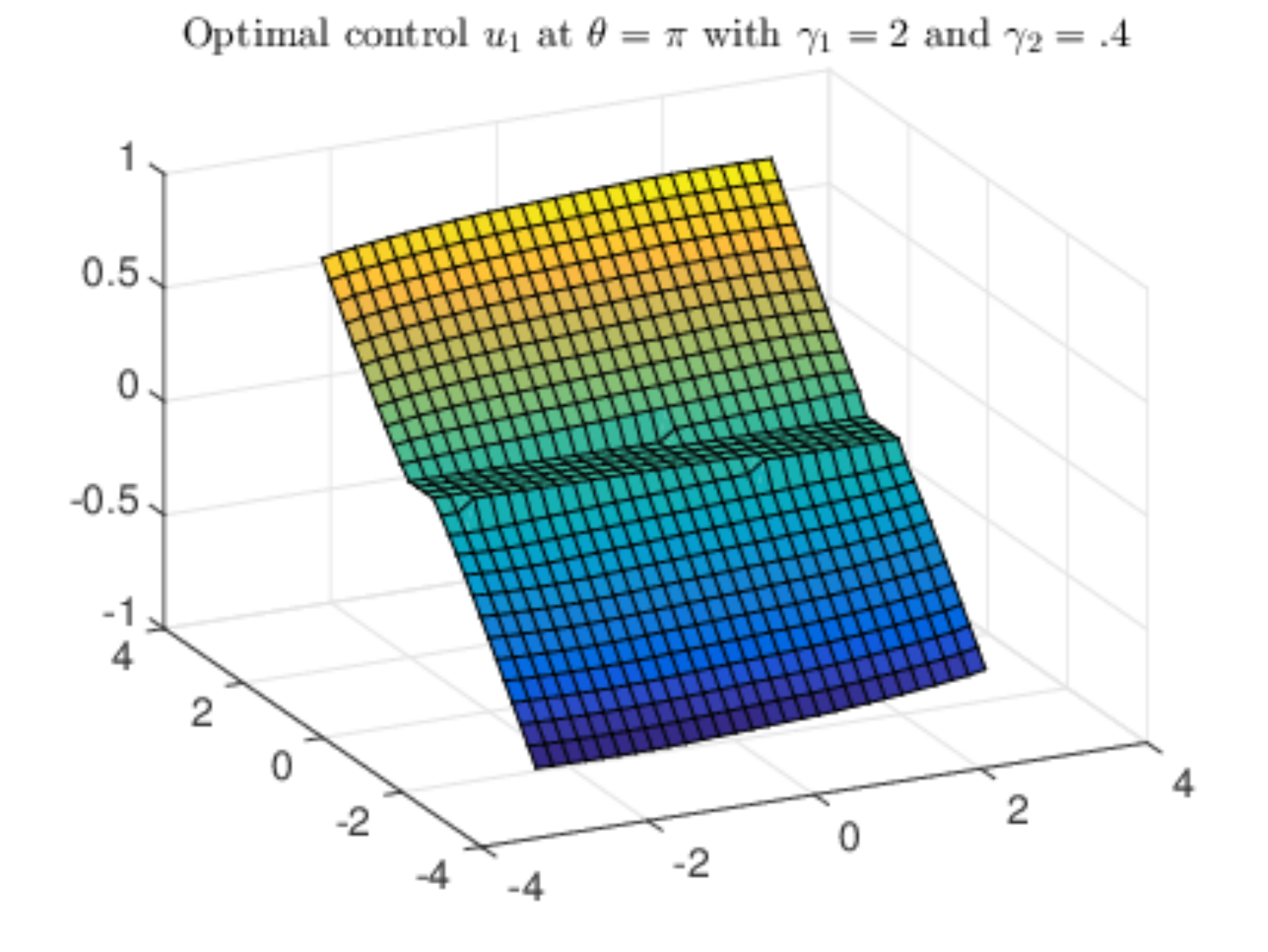,width=.48\linewidth,clip=}
    \epsfig{file=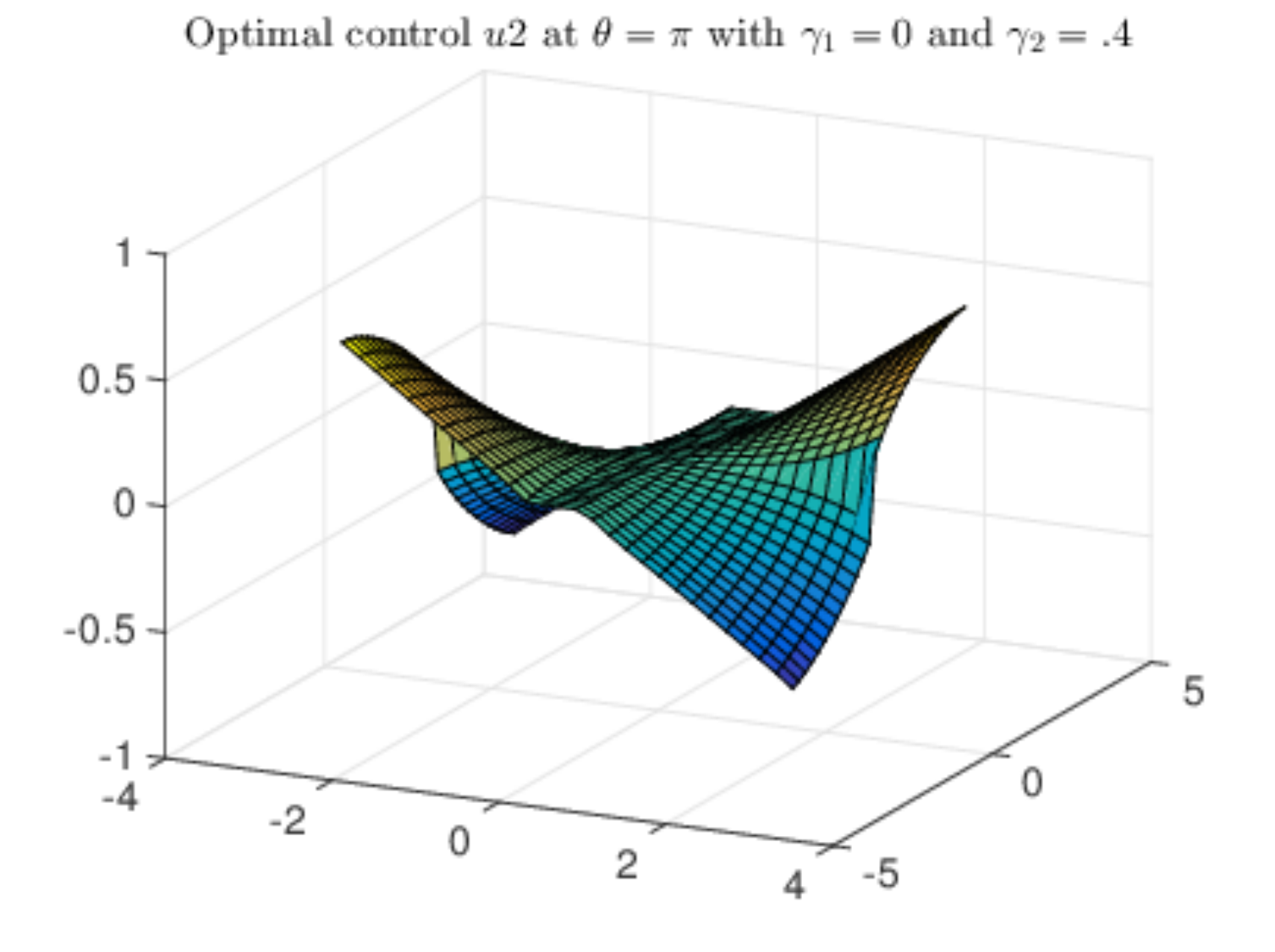,width=.48\linewidth,clip=}
  \epsfig{file=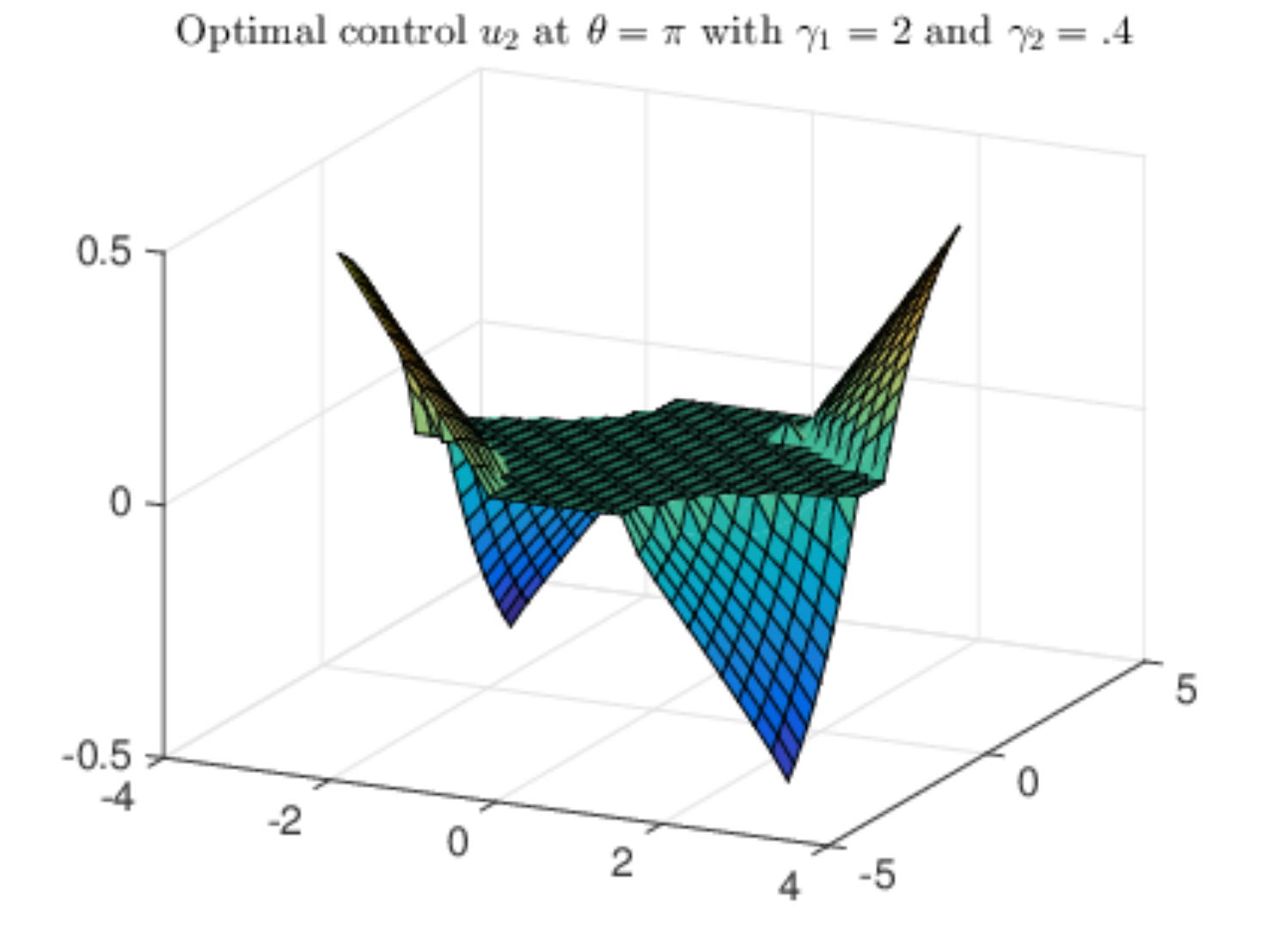,width=.48\linewidth,clip=}\\
    \epsfig{file=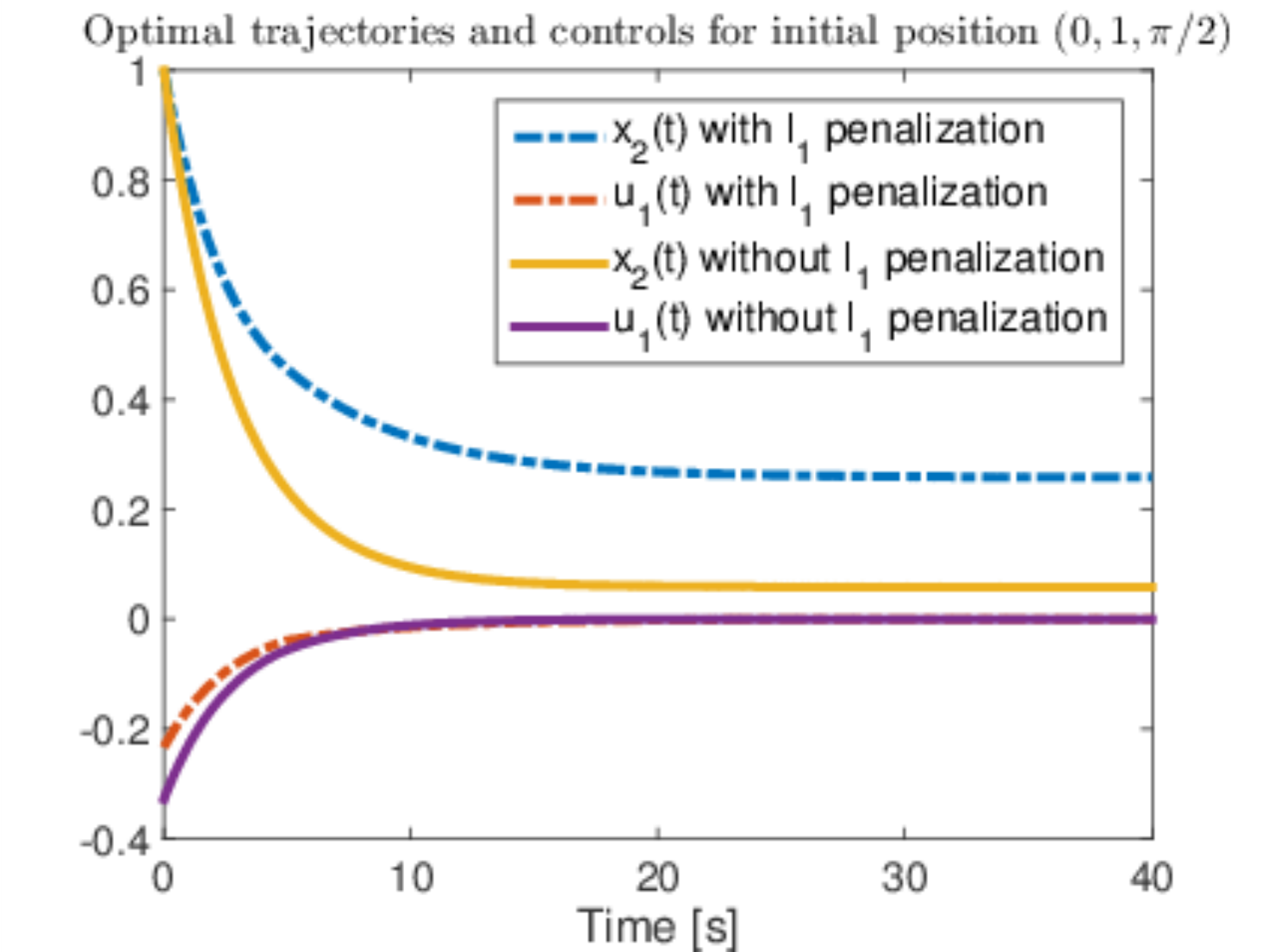,width=.48\linewidth,clip=}\\
  \caption{Sparse control of a 3D car model. Top: different isosurfaces with different $\ell_1$-cost. Rows 2 and 3: effect of sparsity on the optimal control field. Bottom: optimal trajectories with and without sparsity.}\label{test5}
\end{figure}

In our last test, the dynamics are given by a nonlinear 3D car model with two control variables:
\begin{align*}
f(x,u)=(u_1\cos(x_3),u_1\sin(x_3),u_2)
\end{align*}
with $u_1\in\,[-\omega_1\,,\omega_1]$, and $u_2\in\,[-\omega_2,\,\omega_2]$, $\omega_1>0$, $\omega_2 >0$. 
For this problem we implement the semismooth Newton algorithm with box constraints introduced in
Section \ref{sec:ssbc}. Results are shown in Fig. \ref{test5}. The first row shows  same isovalues for different costs. The addition of an $\ell_1$-term shrinks
the region of a given isovalue. The second and third rows depict the effect of sparsity in both control
variables, creating regions of zero action. A direct consequence of this can be seen at the bottom
of Fig. \ref{test5}, where the inclusion of the additional $\ell_1$-cost generates optimal
trajectories which do not reach the origin due to the interplay between the control cost and the discount factor of the optimal control problem.

\section*{Concluding remarks}
We have presented a numerical approach for the solution of HJB equations based on a semi-Lagrangian discretization and the use of different local minimization strategies for the approximation of the corresponding numerical Hamiltonian. The numerical results show a more accurate resolution of the optimal control field at a similar computational cost as the currently used schemes. Furthermore, the proposed approach can be also adapted to treat non-differential costs such as $\ell_1$-penalizations on the control. Since the numerical approximation of the Hamiltonian constitutes one building block within the construction of approximation schemes for HJB and related equations, the idea of using local minimization techniques can be conveniently recast in similar problems, such as front propagation problems and differential games, and in related approximation techniques, like fast marching schemes, policy iteration methods and high-order approximations.

\section*{Appendix A: Decomposition of the control space for Test 3}
We make an explicit presentation of the decomposition of the control space of Test 3. The dynamics
is given by
\[
f(x,u)=(x_2,x_3+u_1,u_2)\,,\quad|u_1|\leq a,\;|u_2|\leq b,\,\quad a,b\in \mathds{R}.
\]
For a given departure point $x^d=(x_1^d,x_2^d,x_3^d)\in\mathds{R}^3$ and a sufficiently small $h$, we want to
identify a relation between subsets $U_{\cI}$ of the control space $U=[-a,a]\times[-b,b]$ and the location of the corresponding arrival points
\[
x^d+h(x_2^d,x_3^d+u_1,u_2)\,,\quad (u_1,u_2)\in U_{\cI}.
\]
In the three-dimensional case, the control set $U$ can be decomposed into at most eight disjoint
subsets, one per octant with
\[
U=\bigcup_{\cI \subset W} U_{\cI}\,
\]
for $W=\{1,2,3 \}$.
Since every octant defines a unique linear interpolant, 
we solve eight different minimization problems, and then compute the global nodal minimizer by comparison. 
For instance, let us consider the sector 
\[Q_{\{1,2,3\}}=\{(x_1,x_2,x_3) \,|\, x_1^d\leq x_1,\,x_2^d\leq x_2,\,x_3^d\leq x_3,\, (x_1-x_1^d)+(x_2-x_2^d)+(x_3-x_3^d)\leq k\}\,,\]
where the evaluation of the arrival point is defined by the linear interpolant $I_{\{1,2,3\}}$ depending on
the points $x^d$,\; $x^d+(k,0,0)$,\; $x^d+(0,k,0)$, and $x^d+(0,0,k)$. The
subset $U_{\{1,2,3\}}$ related to this interpolant reduces to
\[U_{\{1,2,3\}}=\{(u_1,u_2)\in U \,|\, x_3^d+u_1\geq 0,\, u_2\geq 0\}\,.\]
Note that for this definition to be meaningful, we need to assume that $x_2^d\geq 0$, otherwise
$U_{\{1,2,3\}}=\emptyset$. Also note that the condition $(x_1-x_1^d)+(x_2-x_2^d)+(x_3-x_3^d)\leq k$ is omitted by assuming a sufficiently small $h$ ensuring it.
The next subset $U_{\{2,3\}}$ relates to the sector
\[Q_{\{2,3\}}=\{(x_1,x_2,x_3) \,|\, x_1^d\geq x_1,\,x_2^d\leq x_2,\,x_3^d\leq x_3,\, -(x_1-x_1^d)+(x_2-x_2^d)+(x_3-x_3^d)\leq k\}\,,\]
and is given by
 \[
 U_{\{2,3\}}=\{(u_1,u_2)\in U \,|\, x_3^d+u_1\geq 0,\, u_2\geq 0\}\,.
 \] Note that $U_{\{1,2,3\}}\equiv U_{\{2,3\}}$, however, $U_{\{2,3\}}$ is nonempty only when
 $x_2^d\leq 0$. Therefore, for a given departure point, depending on the coordinate $x_2^d$, only one of
 these two subsets will be active with arrival points in different sectors $Q_{\{1,2,3\}}$ or
 $Q_{\{2,3\}}$ (with different interpolation data). In a similar way, the remaining six subsets can be obtained. By assuming a linear control term, the identification of the control subsets is simple, as it will only require the resolution of linear inequalities where the departure point enters as a fixed data.

 \section*{Appendix B: Exact value function for Test 1} The exact value function for the infinite horizon optimal control problem with
eikonal dynamics in Test 1 is derived. The HJB-equation has the form
\begin{equation}\label{hjb}
\lambda v + \underset{u \in U}{\max}\, \left\{ - u^T \nabla v - \left(\frac{\|x\|_2^2}{2} +
\frac{\gamma}{2} \|u\|_2^2\right)\right\}=0,\quad x \in \R^n, \quad \gamma>0,
\end{equation}
where $U = \{u \in \R^2 \, |\, \|u\|_2 \le 1\}, $ for which we want to obtain an explicit solution. If $\frac{1}{\gamma} \|\nabla v\|_2 < 1$, then $u^* = -\frac{1}{\gamma}\nabla v$ provides a maximum for the expression in brackets, and the HJB-equation becomes
\[
\lambda v + \frac{1}{2\gamma}\|\nabla v\|_2^2 - \frac{1}{2} \|x\|_2^2=0.
\]
Switching to polar coordinates $(r, \varphi)$ this equation can be expressed as
\[
\lambda v(r,\varphi) + \frac{1}{2\gamma} (v_r(r,\varphi)^2+ \frac{1}{r^2} v_\varphi(r,\varphi)^2) - \frac{r^2}{2}=0,
\]
and assuming that the solution is radially symmetric
\[
\lambda v(r) + \frac{1}{2\gamma} v_r (r)^2-\frac{r^2}{2} = 0.
\]
The ansatz $v_1=Ar^2$, leads to $\lambda Ar^2 + \frac{2}{\gamma} A^2 r^2 - \frac{r^2}{2} =0$ and
hence $A=\frac{\gamma}{4} (\sqrt{\lambda^2 + \frac{4}{\gamma}}-\lambda).$ The resulting expression
for $v_1$ is the solution of \eqref{hjb}, if $(v_1)_r \le \gamma$ which results in $r \le
2((\lambda^2 + \frac{4}{\gamma})^{\frac{1}{2}}-\lambda)^{-1} = : \overline{r}.$ We note that $v_1(\overline{r})= \frac{\gamma}{2} \overline{r}.$

We turn to the case that $\frac{1}{\gamma} \|\nabla v\|_2 > 1.$ In this case the maximum in \eqref{hjb} is achieved on the boundary of $U$ at $u^* = - \frac{\nabla v}{\|\nabla v\|_2}$. Again we look for a radially symmetric solution in polar coordinates
\begin{equation}\label{2}
\lambda v(r) + |v_r (r)| - \frac{r^2}{2} = \frac{\gamma}{2}.
\end{equation}
Assuming that $v_r \ge 0$ we make the Ansatz $v_2(r) = ar^2 + br + c + de^{-\lambda r}$. Inserting into \eqref{2} and comparing coefficients we obtain
\[
v_2 (r) = \frac{1}{2\lambda} r^2 - \frac{1}{\lambda^2} r + \frac{\gamma}{2\lambda} + \frac{1}{\lambda^3} + d e^{-\lambda r}.
\]
Continuous concatenation with $v_1$ at $\overline{r}$ implies that
\begin{equation}\label{coeff}
d = e^{\lambda \overline{r}} \left[\left( \frac{\gamma}{2} + \frac{1}{\lambda^2}\right) \overline{r} -
\frac{1}{2\lambda} \overline{r}^2 - \frac{\gamma}{2\lambda} - \frac{1}{\lambda^3}\right].
\end{equation}
Note that this latter coupling yields that $v(r)$ is a $\mathcal{C}^1(\R)$ function and therefore it is also a classical solution of eq. \eqref{hjb}.
Summarizing we have
\begin{equation*}
v(r) = \left\{
\begin{aligned}
&\frac{\gamma}{4} (\sqrt{\lambda^2 + \frac{4}{\gamma}} - \lambda)r^2 &&\text{ for } r \le \overline{r}\\
&\frac{1}{2\lambda} r^2 - \frac{1}{\lambda^2} r + \frac{\gamma}{2\lambda} + \frac{1}{\lambda^3} + d e^{-\lambda r} &&\text{ for } r > \overline{r},
\end{aligned}
\right.
\end{equation*}
where $d$ is given in \eqref{coeff}.

\bibliographystyle{amsplain}    
  \bibliography{lit}

\end{document}